\documentclass[12pt,roman]{amsart}
\usepackage{amssymb}
\usepackage{paralist,amsthm,amsfonts,amssymb,graphics, topcapt, amsmath}%
\usepackage{t1enc, pst-poly, pstricks,color,appendix}
\usepackage[active]{srcltx} 
\usepackage{yhmath}

\usepackage{amssymb}
\usepackage{tikz}
\usetikzlibrary{calc}
\usetikzlibrary{decorations.markings}
\usetikzlibrary{patterns,intersections,through}
\usepackage{graphicx}%
\usepackage{yhmath}

\newtheorem{thm}{Theorem}
\newtheorem{lem}{Lemma}
\newtheorem{defn}{Definition}
\newtheorem{cor}{Corollary}

\usepackage{tikz}
\usepackage{graphics}
\usetikzlibrary{calc,arrows}
\usetikzlibrary{decorations.markings}
\usetikzlibrary{patterns,intersections,through}

\begin{document}

\pagestyle{plain}
\pagenumbering{roman}

\pagestyle{plain}
\pagenumbering{roman}

\begin{center}

\vspace*{0.5cm}
\begin{large}
{\bf The Eventual Gaussian Distribution of Self-Intersection Numbers on Closed Surfaces}
\end{large}\\
\vspace*{0.2cm}
by\\
{\large Matthew Wroten}\\
\vspace*{0.2cm}
{\large 2014}
\end{center}
\vspace*{0.4cm}
Abstract:
\\

Oriented loops on an orientable surface are, up to equivalence by free homotopy, in one-to-one correspondence with the conjugacy classes of the surface's fundamental group.  These conjugacy classes can be expressed (not uniquely in the case of closed surfaces) as a cyclic word of minimal length in terms of the fundamental group's generators.  The self-intersection number of a conjugacy class is the minimal number of transverse self-intersections of representatives of the class.  By using Markov chains to encapsulate the exponential mixing of the geodesic flow and achieve sufficient independence, we can use a form of the central limit theorem to describe the statistical nature of the self-intersection number.  For a class chosen at random among all classes of length $n$, the distribution of the self intersection number approaches a Gaussian when n is large.  This theorem generalizes the result of Steven Lalley and Moira Chas to include the case of closed surfaces.

\vspace*{1cm}
\begin{center}
for Amelia Wroten on her twenty-fifth birthday
\end{center}

\tableofcontents

\section*{Acknowledgements}

Foremost I thank God for giving me the ability and life circumstances to complete this work.  Next, I thank Moira Chas and Dennis Sullivan who presented me with this dissertation problem while I was their graduate student at Stony Brook University.  Lastly, thanks also to the Penn State University Mathematics Department for supporting me financially while I verified the positivity of the Gaussian's variance and for Dimitry Dolgopyat for informing me that this only required verifying the function was not cohomologous to a constant.

\pagenumbering{arabic}

\section{Introduction}

In 1987, Marshall Cohen and Martin Lustig studied the paths of closed geodesics on hyperbolic surfaces with boundary \cite{cohen-lustig} \cite{lustig}. They gave a precise combinatorial algorithm for computing the geometric self-intersection number of a closed geodesic. Closed geodesics are in one-to-one correspondence with conjugacy classes of fundamental groups. In 2009, Moira Chas ran a computer program to calculate the self-intersection number of all conjugacy classes of a given word-length for the doubly-punctured plane. She plotted the distribution of the self-intersection numbers of all (approximately) 4,000,000 curves and found it to be remarkably Gaussian in appearance. She consulted Steven Lalley and informed him of the discovery. At first he was doubtful that it could really be a Gaussian until she showed him the histogram. He said that it was the best Gaussian he'd ever seen and so they immediately got to work proving it. In August of 2011, they proved that, for a closed hyperbolic surface with boundary, the distribution of the self-intersection number (when suitably rescaled) approaches a Gaussian as the word-length of the curves approaches infinity \cite{chas-lalley}.

From a combinatorial viewpoint, closed curves are simpler to analyze on a surface with boundary than on a closed surface. This is because for a surface with boundary, the fundamental group is just a free group. On a closed surface, however, there is a relationship among the generators. Because of this, there is not always a unique shortest way to represent a conjugacy class in terms of the generators. In his dissertation in 1987 \cite{lustig}, Martin Lustig described a unique way to codify these conjugacy classes by expressing them in terms of a larger alphabet of $16g$ symbols (as opposed to $4g$ which is the number of generators and inverses for a closed surface of genus $g$) He presented an algorithm to convert a conjugacy class into its unique representative and then calculate the self-intersection number.

It was natural and desirable to do a statistical analysis similar to the one done by Chas and Lalley \cite{chas-lalley} but for closed surfaces. The analogous result does hold here too: for a closed hyperbolic surface, the distribution of self-intersection numbers approaches a Gaussian as the word-length goes to infinity.  First, the conjugacy classes were represented by a Markov process.  Bowen and Series \cite{bowen-series} present a way to do this based on deck transformations.  Here, however, we present a new method that represents more closely the path of the geodesic.  Because of the complicated combinatorial nature of the problem, the Markov chain requires an astonishing $4g(6g-1)$ states. For the simplest case, the surface of genus two, this corresponds to $88$ states. The Parry measure on the Markov chains (the unique measure that maximizes entropy) describes, in a limiting sense, the uniform measure on periodic orbits (conjugacy classes of the fundamental group). Thus by analyzing a statistic such as self-intersection number on the Markov chain, we also know the distribution of the statistic over actual geodesics of a given word-length on the surface. Ergodicity and the exponential mixing of the geodesic flow, which is expressed in terms of the mixing properties of our Markov chain, allow us to achieve a sufficient amount of independence so that a central limit for Markov chains can be utilized.  Finally, we have the limiting normal distribution of the intersection number.

More specifically, let $S$ be a closed orientable surface of genus $g \ge 2$.  Fix a presentatoin of $\pi_1(S)$ such that the universal cover of $S$ can be visualized as a tessalation of the hyperbolic plain by $4g$-gons with $4g$ at each vertex.  Let $N_n(\alpha)$ denote the random variable obtained by evaluating $N$ (the intersection number) on a uniformly chosen $\alpha \in F_n$ (the set of conjugacy classes $\pi_1(S)$ that have (minimal) combinatorial word length $n$ according to the chosen presentation).  Then there exists constants $\kappa$ and $\sigma^2 \ge 0$ (depending only on the genus of $S$ and the presentation of $\pi_1(S)$) such that, as $n \rightarrow \infty$,
$$\frac{N_n(\alpha)-\kappa n^2}{n^{3/2}} \Longrightarrow Normal(0,\sigma)$$
in the weak topology.  This is our main theorem (Theorem \ref{Main}).  Further, Theorem \ref{variance} gives us that $\sigma^2 \gneqq 0$.

\section{The Markov Chain for Snakes}

Throughout the paper (unless otherwise stated) we will explicitly argue the case where the surface $S$ has genus $2$ and its fundamental group is given the presentation $G=<a,b,c,d | abABcdCD=1>$ (here and elsewhere, a capital letter in a group denotes the inverse of the corresponding lowercase letter) , but the arguments generalize to other presentations and genuses as well.  In proving our main theorem, we will need a Markov chain whose periodic orbits correspond to conjugacy classes of the fundamental group of $S$.  Before defining the Markov chain, we will give motivation for its construction by presenting an intuitive classification of the elements of $G$.  Then, after descibing the Markov chain, we will make statements about the way it uniquely (Theorem \ref{uniqueness}) and completely (Theorem \ref{algorithm}) represents conjugacy classes of $G$.  Then we get the required correspondence (Theorem \ref {correspondence}).

\subsection{Snakes and the Fundamental Group}

Our Markov chain will represent geometric limits of what we will call "snakes." Given two vertices in the Cayley graph of the fundamental group of our surface, we define a \emph {snake} to be the union of all shortest paths from one vertex to the other.  The length of a snake is defined as the number of edges in such a shortest path.  We can give snakes a direction by declaring one vertex to be the start and the other vertex to be the finish.  The $n$th step of a snake is defined as the union of the $n$th edges of each of the shortest directed paths from the start to the finish.

Although the the statement of the Main Theorem is purely topological, we will set up a specific geometry (WLOG) to help state and prove some useful theorems.  For now, assume our surface has genus 2 and its fundamental group has the presentation $G=<a,b,c,d | abABcdCD=1>$, and that the surface is formed by identifying appropriate edges of a regular octagon in hyperbolic space with vertex angles $\pi/4$.  We can allow such an octagon, $F$, to be the fundamental domain in the universal cover of the surface.  Tessellating the universal cover with this octagon via the deck transformations of the group yields a graph, $N$ (made of infinitely many infinite geodesics), with geodesic edges in hyperbolic space.  We can embed the Cayley graph, $\Gamma$, of the group dually to $N$ by connecting the centers of all adjacent octagons with geodesics.  Let $v_0$ be the vertex of $\Gamma$ at the center of $F$, corresponding to the identity element of $G$.  There is a canonical bijection between elements of $G$ and vertices of $\Gamma$.  Each vertex in $\Gamma$ also corresponds to the end of a unique snake starting at $v_0$.  Thus the following theorem yields a complete classification of the elements of $G$.

\begin{thm}\label{graph and snakes}
The snakes in $\Gamma$ of length $n$ starting at $v_0$ are in bijective correspondence with the directed paths of length $n+1$ starting at St and finishing at Fi in the directed graph, $G_2$ (see Figure \ref{G_2}).
\end{thm}
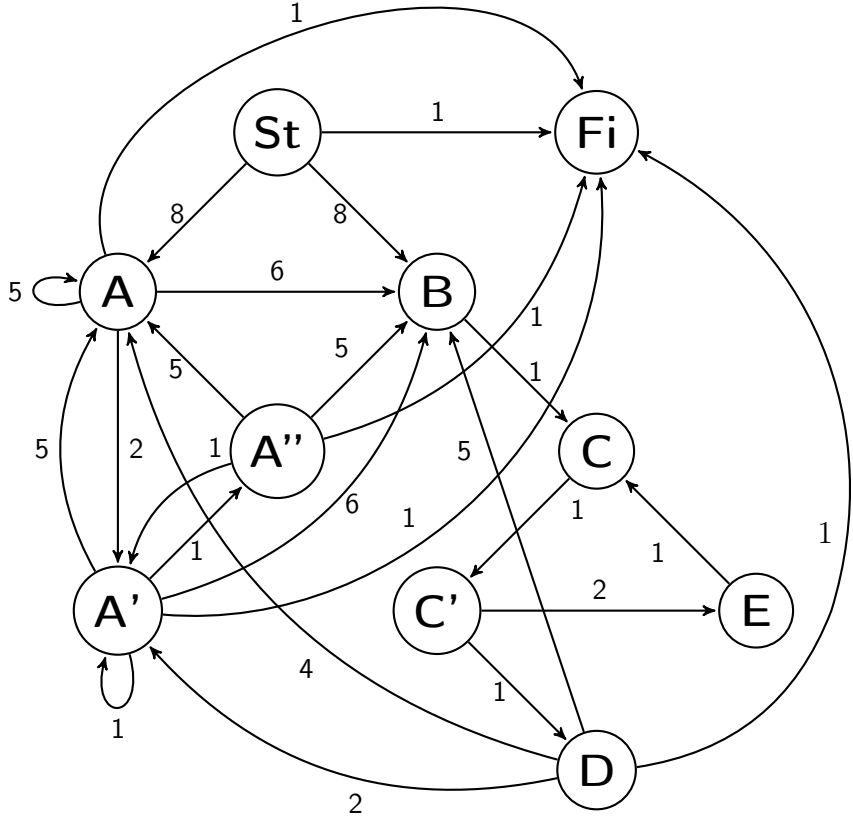
\begin{figure}
\begin{tikzpicture}[->,>=stealth',shorten >=1pt,auto,node distance=3cm,
                    thick,main node/.style={circle,draw,font=\sffamily\Large\bfseries}]

  \node[main node] (1) {St};
  \node[main node] (2) [below left of=1] {A};
  \node[main node] (3) [below right of=1] {B};
  \node[main node] (4) [below right of=2] {A''};
  \node[main node] (5) [below left of=4] {A'};
  \node[main node] (6) [above right of=3] {Fi};
  \node[main node] (7) [below right of=3] {C};
  \node[main node] (8) [below left of=7] {C'};
  \node[main node] (9) [below right of =7] {E};
  \node[main node] (10) [below right of=8] {D};

  \path[every node/.style={font=\sffamily\small}]
    (1) edge node [left] {8} (2)
        edge node[left] {8} (3)
        edge node [above] {1} (6)
    (2) edge node [above] {6} (3)
        edge node {2} (5)
        edge [loop left] node {5} (2)
        edge [bend left = 90] node [above] {1} (6)
    (3) edge node [right] {1} (7)
    (4) edge node [left] {5} (2)
        edge [bend right] node [above, pos = .1] {1} (5)
        edge [bend right] node[right] [pos = .6] {1} (6)
        edge node {5} (3)
     (5) edge [bend left] node {5} (2)
        edge node [right, pos = .3] {1} (4)
        edge [bend right] node [right] {6} (3)
        edge [loop below] node {1} (5)
        edge [bend right = 50] node [pos = .4] {1} (6)
      (7) edge node [pos = .1]{1} (8)
      (8) edge node [left] {1} (10)
         edge node {2} (9)
      (9) edge node {1} (7)
      (10) edge [bend left] node [pos = .4] {4} (2)
         edge [bend left] node [pos = .4] {2} (5)
         edge node [pos = .75] {5} (3);
         \draw (10) to [bend right] ($(9) + (1,0) $ ) to [bend right  = 40] node [pos = .1] {1} (6);
	
\end{tikzpicture}
\caption{This is $G_2$}
\label{G_2}
\end{figure}

Before proving this, we will need a lemma.  Remember that the snake is made up of the union of shortest edgepaths in $\Gamma$ so that knowing the structure of shortest edgepaths will tell us the structure of the snake.  A shortest edgepath between two vertices of $\Gamma$ is a path along the graph of $\Gamma$ from one vertex to the other that traverses the minimal possible number of edges out of all paths in $\Gamma$ between the vertices.  We present a fundamental lemma about shortest edgepaths that will allow us to derive all of their characteristics and thus all of the characteristics of snakes.

\begin{lem}[Golden Lemma]\label{golden lemma}
A shortest edgepath in $\Gamma$ never crosses the same geodesic of $N$ more than once.  Further, a (directed) shortest edgepath can never cross a geodesic $L$ of $N$ into the half-plane bounded by $L$ that doesn't contain its final vertex.
\end{lem}

	Proof:  Suppose $E$ is a shortest edgepath between $v_0$ and $v_n$ that crosses a geodesic, $L$, of $N$ more than once.  Let $v_0, v_1,..., v_n-1, v_n$ be, in order, the vertices traversed by $E$ (we can immediately rule out the possibility that any vertex of $\Gamma$ is traversed multiple times because the portion of the edgepath between the first and last such traversal could be truncated out to yield a shorter edgepath).  Let $v_i$ be the vertex of $E$ immediately preceding the first crossing of $L$.  Let $v_j$ be the vertex of $E$ immediately following the second crossing of $L$.  The portion of $E$ from $v_{i+1}$ to $v_{j-1}$ can be reflected across $L$ to yield a path between $v_i$ and $v_j$ of length $j-i-2$, creating a contradiction.  A (directed) shortest edgepath cannot cross a geodesic $L$ of $N$ into the half-plane bounded by $L$ that doesn't contain its final vertex because it would eventually have to cross back, creating a contradiction $\blacksquare$
	
	Given this lemma, we see that with every edge traversed along a shortest edgepath, we cross a geodesic in $N$, thus eliminating a half-plane from the allowable region that the rest of the edgepath can traverse.  Using this idea, we can find precisely all of the shortest edgepaths from one vertex to another and thus our snake.

	We will now describe the snake that starts at $v_0$ and ends at another vertex $v_n$ of $\Gamma$.  We will break this down into cases based on where $v_n$ lies.  Starting at $v_0$, we will determine, edge by edge, where the shortest edgepath(s) can go based on the Golden Lemma.  For an element of our Markov chain of period n, we want the state of the $i$th step of the element to correspond to the $i$th step in the snake, which is a combination of the $i$th edge in each of the shortest edgepaths together with information about where the next steps can be.  We will build the snake step by step.  At each step we will look at the final vertex or vertices of the edges in the previous step.  We will then decide which edge or edges emanating from there could potentially lead via shortest edgepaths to $v_n$.

Case St: If $v_n$ = $v_0$, we are done (Case Fi) so assume this is not the case.  Then $v_n$ lies outside the fundamental domain $F$.  The 8 geodesics on the boundary of $F$ divide the complement of $F$ into 16 regions -- 8 each of 2 types.

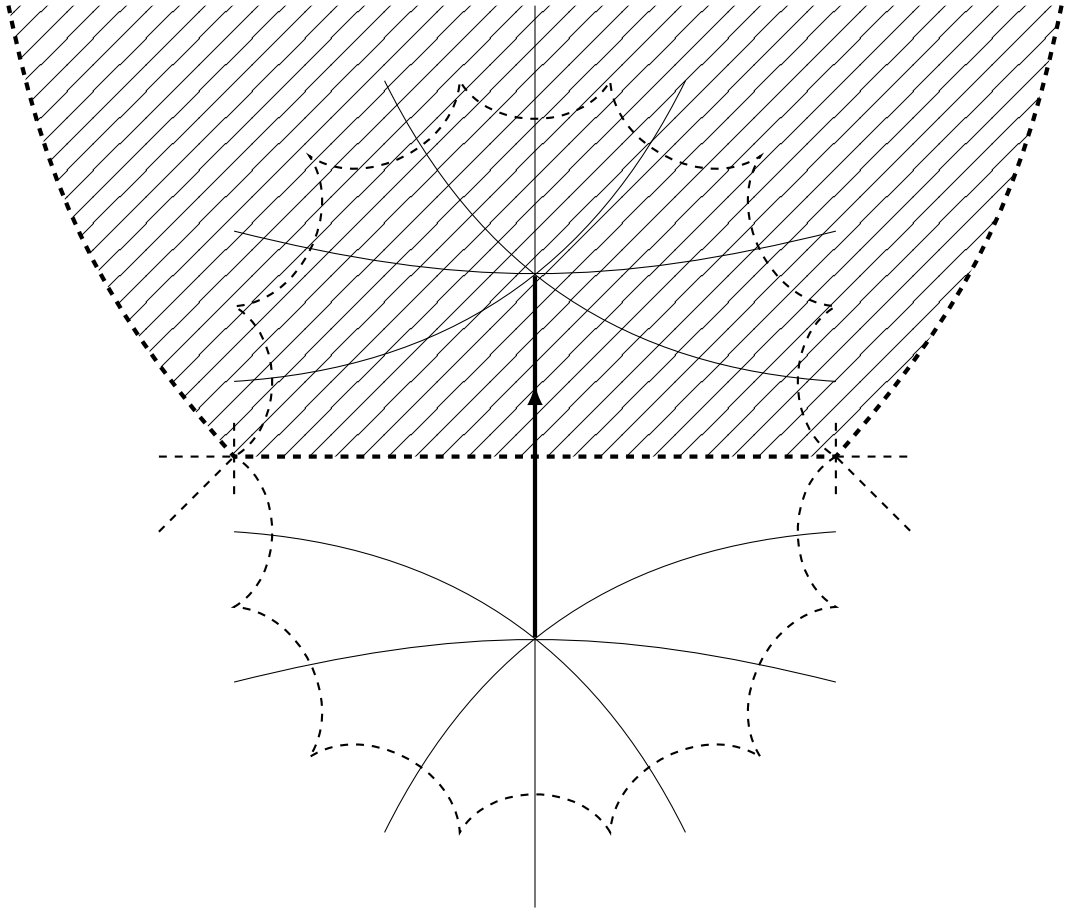
\begin{figure}
\begin{tikzpicture}
\pgfdeclarepatternformonly[\LineSpace]{my north east lines}{\pgfqpoint{-1pt}{-1pt}}{\pgfqpoint{\LineSpace}{\LineSpace}}{\pgfqpoint{\LineSpace}{\LineSpace}}%
{
    \pgfsetlinewidth{0.2pt}
    \pgfpathmoveto{\pgfqpoint{0pt}{0pt}}
    \pgfpathlineto{\pgfqpoint{\LineSpace + 0.1pt}{\LineSpace + 0.1pt}}
    \pgfusepath{stroke}
}

\newdimen\LineSpace
\tikzset{
    line space/.code={\LineSpace=#1},
    line space=10pt
}

\tikzset{
  ->-/.style={
    decoration={
      markings,
      mark=at position #1 with {\arrow[scale=1.5,thick]{latex}},
    },
    postaction={decorate}
  },
  ->-/.default={0.5},
}

	\draw [thick, dashed] (-5,0) -- (-4,0);
	\draw [thick, dashed] (4,0) -- (5,0);
	\draw [thick, dashed, bend right = 60] (4,0) to [bend right] (4,-2) to [bend right] (3,-4) to [bend right] (1,-5) to [bend right] (-1,-5) to [bend right]  (-3, -4) to [bend right] (-4,-2) to [bend right] (-4,0);
	\draw [thick, dashed, bend right = 60] (-5,-1) -- (-4,0) to [bend right] (-4,2) to [bend right] (-3,4) to [bend right] (-1,5) to [bend right] (1,5) to [bend right]  (3, 4) to [bend right] (4,2) to [bend right] (4,0) -- (5,-1);
	\draw [thick, dashed] (4,-.5) -- (4,.5);
	\draw [thick, dashed] (-4,-.5) -- (-4,.5);
  	\draw [ultra thick, dashed, pattern = my north east lines] (7,6) to [bend left = 15] (4,0) -- (-4,0) to [bend left =15] (-7,6);
	\draw (0,-6) -- (0, 6);
	\draw (-4, 1) to [bend right] (2,5);
	\draw (-2,5) to [bend right] (4,1);
	\draw (-4,3) to [bend right = 14] (4,3);
	\draw (4, -1) to [bend right] (-2,-5);
	\draw (2,-5) to [bend right] (-4,-1);
	\draw (4,-3) to [bend right = 14] (-4,-3);
	\draw [ultra thick, ->- = .7] (0,-2.4) -- (0,2.4);

\end{tikzpicture}
\caption{this is a Case A region}
\label{A}
\end{figure}

Case A: $v_n$ is in one of the (open) regions with 3 boundary geodesic pieces, adjacent to $F$.  In this case, any shortest edgepath must immediately cross into that region, otherwise it would cross a different geodesic of $N$ and never be able to enter the region without recrossing it, violating the Golden Lemma (Figure \ref{A}).    We will describe what follows Case A after discussing the other 8 regions.

In all of our diagrams, solid lines are part of $\Gamma$, dashed lines are part of $N$, the edge being discussed will be bold and marked with an arrow, and the region being discussed will be shaded.

Case B: $v_n$ is in one of the (closed) regions with two boundary geodesic rays, adjacent to two of the regions of the Case A type.  In this case, a shortest edgepath can cross first into either of the adjacent regions (but nowhere else) so there can (and will) be multiple edgepaths in the snake (see Figure \ref{B}).  From Case B, the next step in the snake is predetermined.  By the Golden Lemma, the two edgepaths in the snake must continue around the octagon so that they can eventually end up in the Case B region.  Thus we are immediately forced to a Case C.

\begin{figure}
\begin{tikzpicture}[scale =.8]

\newdimen\LineSpace
\tikzset{
    line space/.code={\LineSpace=#1},
    line space=10pt
}

\tikzset{
  ->-/.style={
    decoration={
      markings,
      mark=at position #1 with {\arrow[scale=1.5,thick]{latex}},
    },
    postaction={decorate}
  },
  ->-/.default={0.5},
}

	\draw [ultra thick, dashed, pattern = my north east lines] (8,4) --(0,0) -- (8,-4);
	\draw [thick, dashed, bend right = 60] (3,1.5) to (5,2) to (7,1.5) to (8,0) to (7,-1.5) to (5,-2) to (3,-1.5);
	\draw [thick, dashed] (-3,1.5) to (0,0) to (-3, -1.5);
	\draw [thick, dashed, bend right = 60] (-3,-1.5) to (-5,-2) to (-7,-1.5) to (-8,0) to (-7,1.5) to (-5,2) to (-3,1.5);
	\draw [ultra thick, ->- = .5] (-4,0) to [bend right = 20] (-1.5,2);
	\draw [ultra thick, ->- = .5] (-4,0) to [bend left = 20] (-1.5,-2);
	\draw (-1.5, 2) to [bend right = 60] (0,3) to [bend right = 60] (1.5,2) to [bend right = 20] (4,0) to [bend right = 20] (1.5, -2) to [bend right = 60] (0,-3) to [bend right = 60] (-1.5,-2);
	\draw [thick, dashed] (1,3) to (-1,-3);
	\draw [thick, dashed] (1,-3) to (-1,3); 
	\draw (-8,1) to (-4,0) to (-8,-1);
	\draw (-4,2) to [bend left] (-4,0) to [bend left = 20] (-6.5,-2);
	\draw (-4,-2) to [bend right] (-4,0) to [bend right = 20] (-6.5,2);

\end{tikzpicture}
\caption{this is a Case B region}
\label{B}
\end{figure}

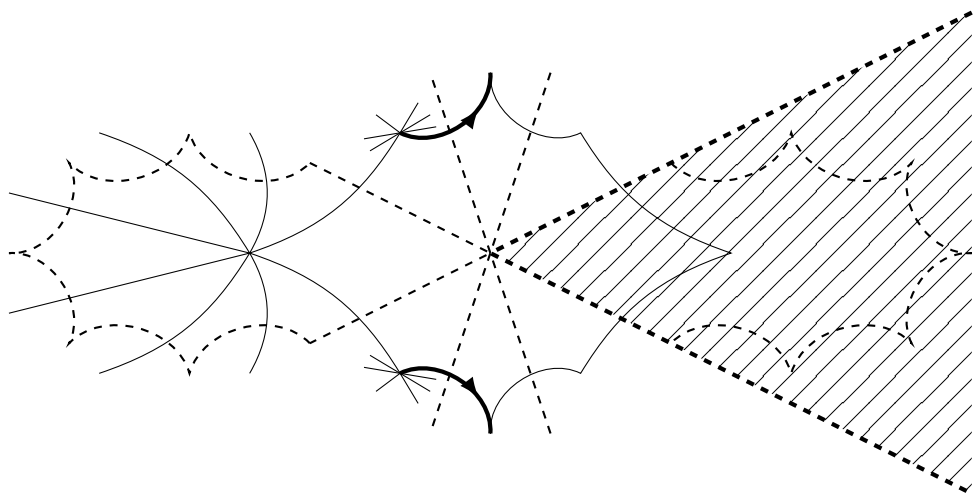
\begin{figure}
\begin{tikzpicture}[scale = .8]

\newdimen\LineSpace
\tikzset{
    line space/.code={\LineSpace=#1},
    line space=10pt
}

\tikzset{
  ->-/.style={
    decoration={
      markings,
      mark=at position #1 with {\arrow[scale=1.5,thick]{latex}},
    },
    postaction={decorate}
  },
  ->-/.default={0.5},
}

	\draw [ultra thick, dashed, pattern = my north east lines] (8,4) --(0,0) -- (8,-4);
	\draw [thick, dashed, bend right = 60] (3,1.5) to (5,2) to (7,1.5) to (8,0) to (7,-1.5) to (5,-2) to (3,-1.5);
	\draw [thick, dashed] (-3,1.5) to (0,0) to (-3, -1.5);
	\draw [thick, dashed, bend right = 60] (-3,-1.5) to (-5,-2) to (-7,-1.5) to (-8,0) to (-7,1.5) to (-5,2) to (-3,1.5);
	\draw [ultra thick, ->- = .7] (-1.5,2) to [bend right = 60] (0,3);
	\draw [ultra thick, ->- = .7] (-1.5,-2) to [bend left = 60] (0,-3);
	\draw (-2,1.7) to (-1,2.3);
	\draw (-2.1,1.9) to (-.9,2.1);
	\draw (-1.9,2.3) to (-1.5,2) to (-1.2, 2.5);
	\draw (-2,-1.7) to (-1,-2.3);
	\draw (-2.1,-1.9) to (-.9,-2.1);
	\draw (-1.9,-2.3) to (-1.5,-2) to (-1.2, -2.5);
	\draw [bend left = 20] (-1.5,2) to (-4,0) to (-1.5,-2);
	\draw (-1.5, 2) to [bend right = 60] (0,3) to [bend right = 60] (1.5,2) to [bend right = 20] (4,0) to [bend right = 20] (1.5, -2) to [bend right = 60] (0,-3) to [bend right = 60] (-1.5,-2);
	\draw [thick, dashed] (1,3) to (-1,-3);
	\draw [thick, dashed] (1,-3) to (-1,3); 
	\draw (-8,1) to (-4,0) to (-8,-1);
	\draw (-4,2) to [bend left] (-4,0) to [bend left = 20] (-6.5,-2);
	\draw (-4,-2) to [bend right] (-4,0) to [bend right = 20] (-6.5,2);
	\

\end{tikzpicture}
\caption{this is a Case C} region
\label{C}
\end{figure}

Case C: $v_n$ is still in the (closed) region described in Case B, but now the edgepaths of the snake need to progress one further step to opposite vertices of an octagon (see Figure \ref{C}).  Again the edgepaths in the snake must continue around the octagon so we are forced into a Case C'.

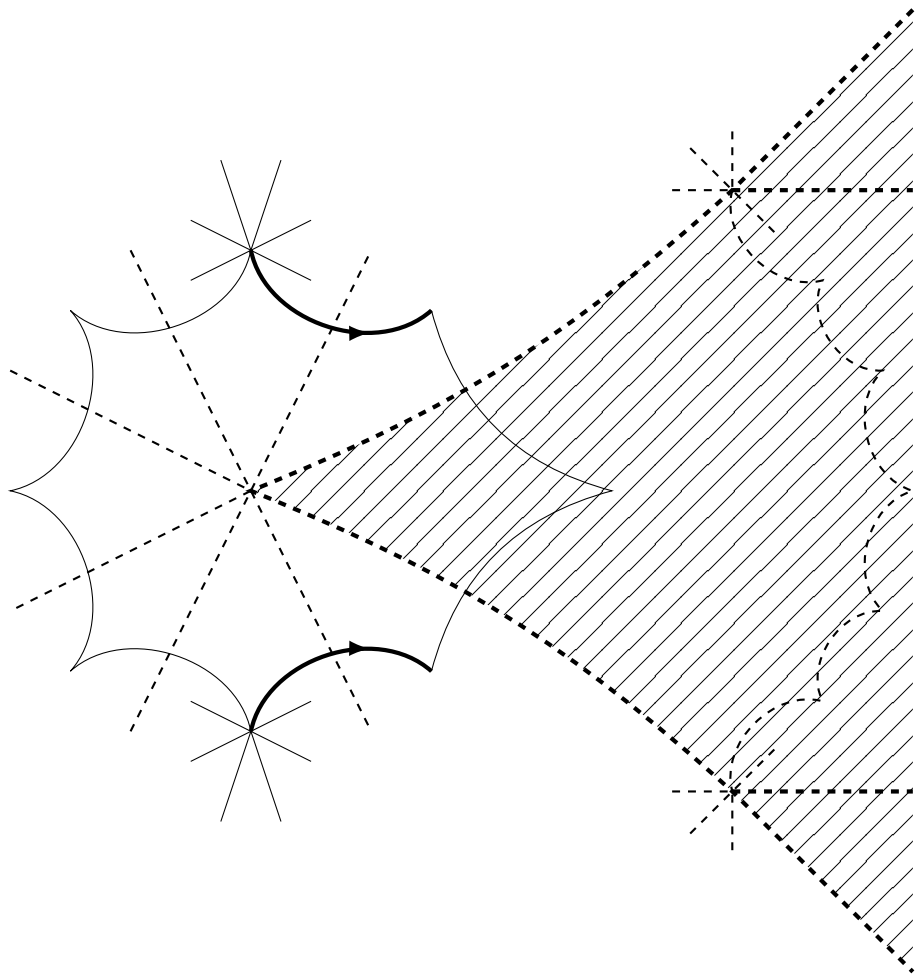
\begin{figure}
\begin{tikzpicture}[scale = .8]

\newdimen\LineSpace
\tikzset{
    line space/.code={\LineSpace=#1},
    line space=10pt
}

\tikzset{
  ->-/.style={
    decoration={
      markings,
      mark=at position #1 with {\arrow[scale=1.5,thick]{latex}},
    },
    postaction={decorate}
  },
  ->-/.default={0.5},
}

	\draw [ultra thick, dashed, pattern = my north east lines] (8,8) to (5,5) to [bend left =10] (-3,0) to [bend left = 10] (5,-5) to (8,-8);
	\draw [ultra thick, dashed] (5,5) to (8,5);
	\draw [ultra thick, dashed] (5,-5) to (8,-5);
	\draw [thick, dashed] (4.3,5.7) to (5.7,4.3);
	\draw [thick, dashed] (4.3,-5.7) to (5.7,-4.3);
	\draw [thick, dashed, bend right = 60] (5,5) to (6.5, 3.5) to (7.5, 2) to (8,0) to (7.5,-2) to (6.5, -3.5) to (5,-5);
	\draw [thick, dashed] (4,5) to (5,5) to (5,6);
	\draw [thick, dashed] (4,-5) to (5,-5) to (5,-6);
	\draw [thick, dashed] (-7,2) to (-3,0) to (-7,-2);
	\draw [thick, dashed] (-5,4) to (-1,-4);
	\draw [thick, dashed] (-5,-4) to(-1,4);
	\draw [ultra thick, ->- = .7] (-3,4) to [bend right = 60] (0,3);
	\draw [ultra thick, ->- = .7] (-3,-4) to [bend left = 60] (0,-3);
	\draw [bend right] (0,3) to (3,0) to (0,-3);
	\draw [bend right = 60] (-3,-4) to (-6,-3) to (-7,0) to (-6,3) to (-3,4);
	\draw (-3.5,5.5) to (-3,4) to (-2.5,5.5);
	\draw (-3.5,-5.5) to (-3, -4) to (-2.5,-5.5);
	\draw (-4, 4.5) to (-2,3.5);
	\draw (-4,3.5) to (-2,4.5);
	\draw (-4, -4.5) to (-2,-3.5);
	\draw (-4,-3.5) to (-2,-4.5);

\end{tikzpicture}
\caption{this is a Case C' region subdivided into its two subregions}
\label{C'}
\end{figure}

Case C': $v_n$ is still in the same (closed) region described as in Case C, but now the edgepaths have traversed one edge closer to (and one edge away from) the final vertex of the aforementioned octagon (by the Golden Lemma they must in order to eventually reach the region) (see  Figure \ref{C'}).  Here there will be choices for the next step depending on where in the region of C' (or C or B) $v_n$ lies.  If $v_n$ is in one of the two wedge shaped regions on the left or right side of the C' region, then we must proceed  to Case E.  Otherwise we proceed to case D.

Case D: $v_n$ is in the (open) region formed by removing two wedges from case E regions.  Here, by the Golden Lemma, both branches of the snake must continue to the final edge of the octagon (see Figure \ref{D}).  Before describing where the snake can go from here, we will continue describing the other cases.

\begin{figure}
\begin{tikzpicture}[scale = .8]

\newdimen\LineSpace
\tikzset{
    line space/.code={\LineSpace=#1},
    line space=10pt
}

\tikzset{
  ->-/.style={
    decoration={
      markings,
      mark=at position #1 with {\arrow[scale=1.5,thick]{latex}},
    },
    postaction={decorate}
  },
  ->-/.default={0.5},
}

	\draw [ultra thick, dashed, pattern = my north east lines] (8,5) to (5,5) to [bend left =10] (-3,0) to [bend left = 10] (5,-5) to (8,-5);
	\draw [ thick, dashed] (5,5) to (6,6);
	\draw [thick, dashed] (5,-5) to (6,-6);
	\draw [thick, dashed] (4.3,5.7) to (5.7,4.3);
	\draw [thick, dashed] (4.3,-5.7) to (5.7,-4.3);
	\draw [thick, dashed, bend right = 60] (5,5) to (6.5, 3.5) to (7.5, 2) to (8,0) to (7.5,-2) to (6.5, -3.5) to (5,-5);
	\draw [thick, dashed] (4,5) to (5,5) to (5,6);
	\draw [thick, dashed] (4,-5) to (5,-5) to (5,-6);
	\draw [thick, dashed] (-7,2) to (-3,0) to (-7,-2);
	\draw [thick, dashed] (-5,4) to (-1,-4);
	\draw [thick, dashed] (-5,-4) to(-1,4);
	\draw [ultra thick, ->- = .7] (0,3) to [bend right] (3,0);
	\draw [ultra thick, ->- = .7] (0,-3) to [bend left] (3,0);
	\draw [bend right] (0,3) to (3,0) to (0,-3);
	\draw [bend right = 60] (0,-3) to (-3,-4) to (-6,-3) to (-7,0) to (-6,3) to (-3,4) to (0,3);
	\draw (-3.5,5.5) to (-3,4) to (-2.5,5.5);
	\draw (-3.5,-5.5) to (-3, -4) to (-2.5,-5.5);
	\draw (-4, 4.5) to (-2,3.5);
	\draw (-4,3.5) to (-2,4.5);
	\draw (-4, -4.5) to (-2,-3.5);
	\draw (-4,-3.5) to (-2,-4.5);
	\draw (-1,3) to (1,3);
	\draw (-.8,3.5) to (.8,2.5);
	\draw (-.3,3.9) to (0,3) to (0.7, 3.5);
	\draw (-1,-3) to (1,-3);
	\draw (-.8,-3.5) to (.8,-2.5);
	\draw (-.3,-3.9) to (0,-3) to (0.7, -3.5);

\end{tikzpicture}
\caption{this is a Case D region}
\label{D}
\end{figure}

Case E:  $v_n$ is now in a (closed) region shaped like the Case B region.  If in the left (right) wedge from Case C', the edgepath on the right (left) branch of the snake must go to the final vertex of the aforementioned octagon.  However, the left (right) branch can get to $v_n$ either via that vertex, or the one to the left (right) (see figure \ref{E}).  From here there will once again be two branches of the snake and we are forced into a Case C.

\begin{figure}
\begin{tikzpicture} [scale = .8]

\newdimen\LineSpace
\tikzset{
    line space/.code={\LineSpace=#1},
    line space=10pt
}

\tikzset{
  ->-/.style={
    decoration={
      markings,
      mark=at position #1 with {\arrow[scale=1.5,thick]{latex}},
    },
    postaction={decorate}
  },
  ->-/.default={0.5},
}

	\draw [thick, dashed] (5,5) to [bend left =10] (-3,0) to [bend left = 10] (5,-5) to (8,-5);
	\draw [ultra thick, dashed, pattern = my north east lines] (8,5) to (5,5) to (8,8);
	\draw [ thick, dashed] (5,5) to (6,6);
	\draw [thick, dashed] (5,-5) to (6,-6);;
	\draw [thick, dashed] (3,7) to (5, 5) to [bend right] (7.5,4);
	\draw [thick, dashed, bend right = 60] (5,5) to (6.5, 3.5) to (7.5, 2) to (8,0) to (7.5,-2) to (6.5, -3.5) to (5,-5);
	\draw [thick, dashed] (2,5) to (5,5) to (5,8);
	\draw [thick, dashed] (4,-5) to (5,-5) to (5,-6);
	\draw [thick, dashed] (-7,2) to (-3,0) to (-7,-2);
	\draw [thick, dashed] (-5,4) to (-1,-4);
	\draw [thick, dashed] (-5,-4) to(-1,4);
	\draw [ultra thick, ->- = .7] (0,3) to (3,0);
	\draw [ultra thick, ->- = .7] (0,-3) to [bend left] (3,0);
	\draw [bend right = 60] (0,-3) to (-3,-4) to (-6,-3) to (-7,0) to (-6,3) to (-3,4) to (0,3);
	\draw [ultra thick, ->- =.5 ] (0,3) to [ bend right=70] (3,6);
	\draw (-3.5,5.5) to (-3,4) to (-2.5,5.5);
	\draw (-3.5,-5.5) to (-3, -4) to (-2.5,-5.5);
	\draw (-4, 4.5) to (-2,3.5);
	\draw (-4,3.5) to (-2,4.5);
	\draw (-4, -4.5) to (-2,-3.5);
	\draw (-4,-3.5) to (-2,-4.5);
	\draw [bend right] (3,6) to (4,7) to (6,7) to (7,6) to (7,4.5) to (6,3.8) to (3,0);
	\draw (-1,-3) to (1,-3);
	\draw (-.8,-3.5) to (.8,-2.5);
	\draw (-.3,-3.9) to (0,-3) to (0.7, -3.5);
	\draw (-1,2.8) to (1,3.2);
	\draw (-.7,3.7) to (0,3) to (.7,3.7);
	\draw (0,3) to (-1,3.2);

\end{tikzpicture}
\caption{this is a Case E region}
\label{E}
\end{figure}

	We return to the discussion of Case A.  Once in case A, there are 7 edges that can be traversed during the next step (the eighth edge (the one we just traversed) is ruled out by the Golden Lemma).  The edge(s) taken next is(are) determined by where in the region of Case A $v_n$ lies.  We consider the octagon of $N$ containing the final vertex of the edge just taken.  If we extend the geodesic segments that make up the boundary of this octagon, the region of Case A is divided into 14 pieces (see Figure \ref{A cont}).  One is the current octagon.  If $v_n$ is there, then we are done (Case Fi).  6 of the regions correspond to entering Case B, and 5 of the regions correspond to entering another Case A.  The regions on the left and right extreme are almost like Case A, but a wedge has been removed from the far side (when in the previous step we assumed we were in Case A and not B).  We will call such a region Case A'.

\begin{figure}
\begin{tikzpicture}[scale = .7]

\newdimen\LineSpace
\tikzset{
    line space/.code={\LineSpace=#1},
    line space=10pt
}

\tikzset{
  ->-/.style={
    decoration={
      markings,
      mark=at position #1 with {\arrow[scale=1.5,thick]{latex}},
    },
    postaction={decorate}
  },
  ->-/.default={0.5},
}

	\draw [thick, dashed] (-5,0) -- (-4,0);
	\draw [thick, dashed] (4,0) -- (5,0);
	\draw [thick, dashed, bend right = 60] (4,0) to [bend right] (4,-2) to [bend right] (3,-4) to [bend right] (1,-5) to [bend right] (-1,-5) to [bend right]  (-3, -4) to [bend right] (-4,-2) to [bend right] (-4,0);
	\draw [ultra thick, dashed, bend right = 60] (-4,0) to [bend right] (-4,2) to [bend right] (-3,4) to [bend right] (-1,5) to [bend right] (1,5) to [bend right]  (3, 4) to [bend right] (4,2) to [bend right] (4,0);
	\draw [ultra thick, dashed] (-8,6) -- (-4,2) to [bend left = 20] (-9,6);
	\draw [ultra thick, dashed] (-4,6) to [bend left = 10] (-3,4) to [bend left =10] (-5,6);
	\draw [ultra thick, dashed] (-1,6) -- (-1,5) -- (-2,6);
	\draw [thick, dashed] (-.5, 5.5) -- (-1.5,4.5);
	\draw [thick, dashed] (-1.5, 5) -- (-.5,5);
	\draw [thick, dashed] (-3.5,3.8) -- (-2.5, 4.2);
	\draw [thick, dashed] (-3.2,3.5) -- (-2.8, 4.5);
	\draw [thick, dashed] (-4,2.5) -- (-4,1.5);
	\draw [thick, dashed] (-3.5,2.5)--(-4.5,1.5);
	\draw [ultra thick, dashed] (8,6) -- (4,2) to [bend right = 20] (9,6);
	\draw [ultra thick, dashed] (4,6) to [bend right = 10] (3,4) to [bend right =10] (5,6);
	\draw [ultra thick, dashed] (1,6) -- (1,5) -- (2,6);
	\draw [thick, dashed] (.5, 5.5) -- (1.5,4.5);
	\draw [thick, dashed] (1.5, 5) -- (.5,5);
	\draw [thick, dashed] (3.5,3.8) -- (2.5, 4.2);
	\draw [thick, dashed] (3.2,3.5) -- (2.8, 4.5);
	\draw [thick, dashed] (4,2.5) -- (4,1.5);
	\draw [thick, dashed] (3.5,2.5)--(4.5,1.5);
	\draw [thick, dashed] (4,-.5) -- (4,.5);
	\draw [thick, dashed] (-4,-.5) -- (-4,.5);
  	\draw [ultra thick, dashed, pattern = my north east lines] (10,6) to [bend left = 15] (4,0) -- (-4,0) to [bend left =15] (-10,6);
	\draw (0,-6) -- (0, 6);
	\draw (-4, 1) to [bend right] (2,5);
	\draw (-2,5) to [bend right] (4,1);
	\draw (-4,3) to [bend right = 14] (4,3);
	\draw (4, -1) to [bend right] (-2,-5);
	\draw (2,-5) to [bend right] (-4,-1);
	\draw (4,-3) to [bend right = 14] (-4,-3);
	\draw [ultra thick, ->- = .7] (0,-2.4) -- (0,2.4);
	\draw (-5,4) node {A};
	\draw (5,4) node {A};
	\draw (-2.5,5) node {A};
	\draw (2.5,5) node {A};
	\draw (0,5.2) node {A};
	\draw (-5,1.5) node {A'};
	\draw (5,1.5) node {A'};
	\draw (-6.8, 4) node {B};
	\draw (6.8,4) node {B};
	\draw (-3.8,5.2) node {B};
	\draw (3.8,5.2) node {B};
	\draw (-1.3,5.7) node {B};
	\draw (1.3, 5.7) node{B};

\end{tikzpicture}
\caption{this is a Case A region subdivided into 14 pieces corresponding to the different possible subsequent regions}
\label{A cont}
\end{figure}

Case A':  $v_n$ is in a (open) region shaped like a Case A region but with a wedge removed from the left (right) (see Figure \ref{A'}).  By the Golden Lemma, the snake must immediately cross into this region.  Again we extend the geodesics of the octagon of N surrounding the final vertex of the edge just taken.  This divides the A' region again into 14 pieces.  If $v_n$ is in the octagon, we are done (Case Fi).  6 regions lead to Case B and 5 regions lead to case A.  The region on the extreme right (left) leads to another Case A', but the region on the extreme left (right) will lead to a narrower A''.

\begin{figure}
\begin{tikzpicture}[scale = .7]

\newdimen\LineSpace
\tikzset{
    line space/.code={\LineSpace=#1},
    line space=10pt
}

\tikzset{
  ->-/.style={
    decoration={
      markings,
      mark=at position #1 with {\arrow[scale=1.5,thick]{latex}},
    },
    postaction={decorate}
  },
  ->-/.default={0.5},
}

	\draw [thick, dashed] (-5,0) -- (-4,0);
	\draw [thick, dashed] (4,0) -- (5,0);
	\draw [thick, dashed, bend right = 60] (4,0) to [bend right] (4,-2) to [bend right] (3,-4) to [bend right] (1,-5) to [bend right] (-1,-5) to [bend right]  (-3, -4) to [bend right] (-4,-2) to [bend right] (-4,0);
	\draw [ultra thick, dashed, bend right = 60] (-4,0) to [bend right] (-4,2) to [bend right] (-3,4) to [bend right] (-1,5) to [bend right] (1,5) to [bend right]  (3, 4) to [bend right] (4,2) to [bend right] (4,0);
	\draw [ultra thick, dashed] (-8,6) -- (-4,2) to [bend left = 20] (-9,6);
	\draw [ultra thick, dashed] (-4,6) to [bend left = 10] (-3,4) to [bend left =10] (-5,6);
	\draw [ultra thick, dashed] (-1,6) -- (-1,5) -- (-2,6);
	\draw [thick, dashed] (-.5, 5.5) -- (-1.5,4.5);
	\draw [thick, dashed] (-1.5, 5) -- (-.5,5);
	\draw [thick, dashed] (-3.5,3.8) -- (-2.5, 4.2);
	\draw [thick, dashed] (-3.2,3.5) -- (-2.8, 4.5);
	\draw [thick, dashed] (-4,2.5) -- (-4,1.5);
	\draw [thick, dashed] (-3.5,2.5)--(-4.5,1.5);
	\draw [ultra thick, dashed] (8,6) -- (4,2) to [bend right = 20] (9,6);
	\draw [ultra thick, dashed] (4,6) to [bend right = 10] (3,4) to [bend right =10] (5,6);
	\draw [ultra thick, dashed] (1,6) -- (1,5) -- (2,6);
	\draw [thick, dashed] (.5, 5.5) -- (1.5,4.5);
	\draw [thick, dashed] (1.5, 5) -- (.5,5);
	\draw [thick, dashed] (3.5,3.8) -- (2.5, 4.2);
	\draw [thick, dashed] (3.2,3.5) -- (2.8, 4.5);
	\draw [thick, dashed] (4,2.5) -- (4,1.5);
	\draw [thick, dashed] (3.5,2.5)--(4.5,1.5);
	\draw [thick, dashed] (4,-.5) -- (4,.5);
	\draw [thick, dashed] (-4,-.5) -- (-4,0);
	\draw [thick, dashed] (-4.5, .5) -- (-4,0) -- (-4.5,-.5);
	\draw [thick, dashed] (4,0) -- (4.5,-.5);
  	\draw [ultra thick, dashed, pattern = my north east lines] (10,6) to [bend left = 15] (4,0) -- (-4,0) to [bend right =40] (-4.7,1) to [bend left = 15] (-9.7,6);
	\draw (0,-6) -- (0, 6);
	\draw (-4, 1) to [bend right] (2,5);
	\draw (-2,5) to [bend right] (4,1);
	\draw (-4,3) to [bend right = 14] (4,3);
	\draw (4, -1) to [bend right] (-2,-5);
	\draw (2,-5) to [bend right] (-4,-1);
	\draw (4,-3) to [bend right = 14] (-4,-3);
	\draw [ultra thick, ->- = .7] (0,-2.4) -- (0,2.4);
	\draw (-5,4) node {A};
	\draw (5,4) node {A};
	\draw (-2.5,5) node {A};
	\draw (2.5,5) node {A};
	\draw (0,5.2) node {A};
	\draw (-5,1.8) node {A''};
	\draw (5,1.5) node {A'};
	\draw (-6.8, 4) node {B};
	\draw (6.8,4) node {B};
	\draw (-3.8,5.2) node {B};
	\draw (3.8,5.2) node {B};
	\draw (-1.3,5.7) node {B};
	\draw (1.3, 5.7) node{B};

\end{tikzpicture}
\caption{this is a Case A' region subdivided into 14 pieces corresponding to the different possible subsequent regions}
\label{A'}
\end{figure}

Case A'': This is (open) shaped like Case A' but with the final edge on the left (right) removed.  Again by the Golden Lemma, the snake must cross into this region immediately.  Dividing the region as before yields only 12 regions (see Figure \ref{A''}).  If $v_n$ is in the octagon, we're done (Case Fi).  Apart from that, 5 regions (including the one on the far left (right)) lead to Case A, 5 regions lead to Case B, and the region on the far right (left) leads to another Case A'.

\begin{figure}
\begin{tikzpicture}[scale = .7]

\newdimen\LineSpace
\tikzset{
    line space/.code={\LineSpace=#1},
    line space=10pt
}

\tikzset{
  ->-/.style={
    decoration={
      markings,
      mark=at position #1 with {\arrow[scale=1.5,thick]{latex}},
    },
    postaction={decorate}
  },
  ->-/.default={0.5},
}

	\draw [thick, dashed] (-5,0) -- (-4,0);
	\draw [thick, dashed] (4,0) -- (5,0);
	\draw [thick, dashed, bend right = 60] (4,0) to [bend right] (4,-2) to [bend right] (3,-4) to [bend right] (1,-5) to [bend right] (-1,-5) to [bend right]  (-3, -4) to [bend right] (-4,-2) to [bend right] (-4,0);
	\draw [ultra thick, dashed, bend right = 60] (-4,2) to [bend right] (-3,4) to [bend right] (-1,5) to [bend right] (1,5) to [bend right]  (3, 4) to [bend right] (4,2) to [bend right] (4,0);
	\draw [ultra thick, dashed] (-4,6) to [bend left = 10] (-3,4) to [bend left =10] (-5,6);
	\draw [ultra thick, dashed] (-1,6) -- (-1,5) -- (-2,6);
	\draw [thick, dashed] (-.5, 5.5) -- (-1.5,4.5);
	\draw [thick, dashed] (-1.5, 5) -- (-.5,5);
	\draw [thick, dashed] (-3.5,3.8) -- (-2.5, 4.2);
	\draw [thick, dashed] (-3.2,3.5) -- (-2.8, 4.5);
	\draw [thick, dashed] (-4,2.5) -- (-4,1.5);
	\draw [thick, dashed] (-3.5,2.5)--(-4.5,1.5);
	\draw [ultra thick, dashed] (8,6) -- (4,2) to [bend right = 20] (9,6);
	\draw [ultra thick, dashed] (4,6) to [bend right = 10] (3,4) to [bend right =10] (5,6);
	\draw [ultra thick, dashed] (1,6) -- (1,5) -- (2,6);
	\draw [thick, dashed] (.5, 5.5) -- (1.5,4.5);
	\draw [thick, dashed] (1.5, 5) -- (.5,5);
	\draw [thick, dashed] (3.5,3.8) -- (2.5, 4.2);
	\draw [thick, dashed] (3.2,3.5) -- (2.8, 4.5);
	\draw [thick, dashed] (4,2.5) -- (4,1.5);
	\draw [thick, dashed] (3.5,2.5)--(4.5,1.5);
	\draw [thick, dashed] (4,-.5) -- (4,.5);
	\draw [thick, dashed] (-4,-.5) -- (-4,0);
	\draw [thick, dashed] (-4.5, .5) -- (-4,0) -- (-4.5,-.5);
	\draw [thick, dashed] (4,0) -- (4.5,-.5);
  	\draw [ultra thick, dashed, pattern = my north east lines] (10,6) to [bend left = 15] (4,0) -- (-4,0) to [bend right=60] (-4,2) to (-8,6);
	\draw (0,-6) -- (0, 6);
	\draw (-4, 1) to [bend right] (2,5);
	\draw (-2,5) to [bend right] (4,1);
	\draw (-4,3) to [bend right = 14] (4,3);
	\draw (4, -1) to [bend right] (-2,-5);
	\draw (2,-5) to [bend right] (-4,-1);
	\draw (4,-3) to [bend right = 14] (-4,-3);
	\draw [thick, dashed] (-4,0) -- (-4, .5);
	\draw [thick, dashed] (-4,2) -- (-4.5,2);
	\draw [ultra thick, ->- = .7] (0,-2.4) -- (0,2.4);
	\draw (-5,4) node {A};
	\draw (5,4) node {A};
	\draw (-2.5,5) node {A};
	\draw (2.5,5) node {A};
	\draw (0,5.2) node {A};
	\draw (5,1.5) node {A'};
	\draw (6.8,4) node {B};
	\draw (-3.8,5.2) node {B};
	\draw (3.8,5.2) node {B};
	\draw (-1.3,5.7) node {B};
	\draw (1.3, 5.7) node{B};

\end{tikzpicture}
\caption{this is a Case A'' region subdivided into 12 pieces corresponding to the different possible subsequent regions}
\label{A''}
\end{figure}

	We conclude by finishing the discussion of Case D.  After splitting up the region by extending geodesics segments of the boundary of the octagon of N containing the vertex where the previous to branches of the snake merged, there are 12 pieces (see Figure \ref{D cont}).  If $v_n$ is in the octagon in the middle, we're done (Case Fi).  Apart from that, 5 regions lead to Case B, 4 regions lead to Case A, and the regions on the extreme left and right lead to Case A'.

\begin{figure}
\begin{tikzpicture}[scale =.8]

\newdimen\LineSpace
\tikzset{
    line space/.code={\LineSpace=#1},
    line space=10pt
}

\tikzset{
  ->-/.style={
    decoration={
      markings,
      mark=at position #1 with {\arrow[scale=1.5,thick]{latex}},
    },
    postaction={decorate}
  },
  ->-/.default={0.5},
}

	\draw [ultra thick, dashed, pattern = my north east lines] (9,7) to [bend left = 30] (5,5) to [bend left =10] (-3,0) to [bend left = 10] (5,-5) to [bend left = 30] (9,-7);
	\draw [ thick, dashed] (5,5) to (6,6);
	\draw [thick, dashed] (5,-5) to (6,-6);
	\draw [thick, dashed] (4.3,5.7) to (5.5,4.5);
	\draw [thick, dashed] (4.3,-5.7) to (5.5,-4.5);
	\draw [ultra thick, dashed, bend right = 50] (5,5) to (6.5, 3.5) to (7.5, 2) to (8,0) to (7.5,-2) to (6.5, -3.5) to (5,-5);
	\draw [ultra thick, dashed] (9,3.5) -- (6.5,3.5) to [bend left =20] (9,5);
	\draw [ultra thick, dashed] (9, 1.8) --(7.5,2) to [bend left = 20] (9,2.8);
	\draw [ultra thick, dashed] (9, .6) -- (8,0) -- (9,-.6);
	\draw [ultra thick, dashed] (9,-3.5) -- (6.5,-3.5) to [bend right =20] (9,-5);
	\draw [ultra thick, dashed] (9, -1.8) --(7.5,-2) to [bend right = 20] (9,-2.8);
	\draw [thick, dashed] (4,5) to (5,5) to (5,6);
	\draw [thick, dashed] (4,-5) to (5,-5) to (5,-6);
	\draw [thick, dashed] (-7,2) to (-3,0) to (-7,-2);
	\draw [thick, dashed] (-5,4) to (-1,-4);
	\draw [thick, dashed] (-5,-4) to(-1,4);
	\draw [ultra thick, ->- = .7] (0,3) to [bend right] (3,0);
	\draw [ultra thick, ->- = .7] (0,-3) to [bend left] (3,0);
	\draw [bend right] (0,3) to (3,0) to (0,-3);
	\draw [bend right = 60] (0,-3) to (-3,-4) to (-6,-3) to (-7,0) to (-6,3) to (-3,4) to (0,3);
	\draw (8,1) -- (3,0) -- (8,-1);
	\draw (7.5, 3) to [bend right] (3,0) to [bend right] (7.5,-3);
	\draw (5.7,4.3) to [bend right=50] (3,0) to [bend right=50] (5.7,-4.3);
	\draw (-3.5,5.5) to (-3,4) to (-2.5,5.5);
	\draw (-3.5,-5.5) to (-3, -4) to (-2.5,-5.5);
	\draw (-4, 4.5) to (-2,3.5);
	\draw (-4,3.5) to (-2,4.5);
	\draw (-4, -4.5) to (-2,-3.5);
	\draw (-4,-3.5) to (-2,-4.5);
	\draw (-1,3) to (1,3);
	\draw (-.8,3.5) to (.8,2.5);
	\draw (-.3,3.9) to (0,3) to (0.7, 3.5);
	\draw (-1,-3) to (1,-3);
	\draw (-.8,-3.5) to (.8,-2.5);
	\draw (-.3,-3.9) to (0,-3) to (0.7, -3.5);
	\draw [thick, dashed] (6.3, 4) -- (6.7,3);
	\draw [thick, dashed] (6,3.8) -- (7, 3.2);
	\draw [thick, dashed] (7.5,2.5) -- (7.5,1.5);
	\draw [thick, dashed] (7.2, 2.4) -- (7.8,1.6);
	\draw [thick, dashed] (7.8,.5) -- (8.2,-.5);
	\draw [thick, dashed] (7.8,-.5) -- (8.2, .5);
	\draw [thick, dashed] (6.3, -4) -- (6.7,-3);
	\draw [thick, dashed] (6,-3.8) -- (7, -3.2);
	\draw [thick, dashed] (7.5,-2.5) -- (7.5,-1.5);
	\draw [thick, dashed] (7.2, -2.4) -- (7.8,-1.6);
	\draw (8.5, 5.5) node {A'};
	\draw (8.5,4) node {B};
	\draw (8,3) node {A};
	\draw (8.5, 2.2) node {B};
	\draw (8.2, 1) node {A};
	\draw (8.6, 0) node {B};
	\draw (8.5, -5.5) node {A'};
	\draw (8.5,-4) node {B};
	\draw (8,-3) node {A};
	\draw (8.5, -2.2) node {B};
	\draw (8.2, -1) node {A};

\end{tikzpicture}
\caption{this is a Case D region subdivided into 12 pieces corresponding to the different possible subsequent regions}
\label{D cont}
\end{figure}

	$G_2$ is constructed from the case by case scenario above as follows.  Each Case is given a vertex of the graph with the same label.  If there are $x$ ways that Case $Y$ leads directly to Case $Z$ (corresponding to $x$ Case $Z$ shaped regions in a Case $Y$ shaped region), then we draw $x$ directed edges from vertex $Y$ to vertex $Z$.  We label each of these edges $YZ_1, YZ_2,... ,YZ_x$.  We say $YZ_1$ represents entering the Case $Z$ shaped region furthest to the left relative to the direction the snake is progressing (Excepting with the StA and StB edges where an arbitrary starting point must be chosen).  Then the rest of $YZ_2,..., YZ_x$ represent in clockwise order the other Case $Z$ shaped regions.  

	Proof of Theorem \ref{graph and snakes}:  First we show how to construct the snake corresponding to a directed path $P$ in $G_2$ from St to Fi.  Start at $v_0$.  Whenever we add an edge to $P$, this will correspond to the addition of a step to the snake and the entering of the corresponding region, unless the edge leads to Fi, in which case the snake stops there (we immediately have the correspondence between the length $n+1$ of a path in $G_2$ and the length n of a snake).  Each successive region in nonempty, and when Fi is reached, a specific vertex is chosen in $\Gamma$.  Thus our correspondence is well-defined.  To check that the correspondence is injective, one only needs to observe that in each Case, the region is split into non-overlapping sub-regions.  All that remains is to check whether the correspondence is onto.

	Choose a vertex $v$ in $\Gamma$.  $\Gamma$ is connected so there exist finite edgepaths from $v_0$ to $v$.  All of these edgepaths have finite, non-negative, integer lengths and since the non-negative integers are well ordered, there exists a (not necessarily unique) edgepath of minimal length.  Call this length $n$.  Proceed through the case by case procedure above, building a snake from $v_0$ toward $v$.  This cannot end before the $n+1$ step or we would have a path of length less than $n$.  At step $j$, our procedure has constructed all edgepaths of length $j$ that could possibly be a part of a shortest edgepath from $v_0$ to $v$ by pruning only edges (and thus paths) that are ruled out by the Golden Lemma.  Thus, at least one of the end-vertices of the $n$th step of the snake must be $v$.  Suppose that there are two end-vertices of the $n$th step of the snake.  Then the $n$th step would have to be a Case B, Case C, Case C', or Case E.  However in each of these cases, it is assumed that the end-vertices of the step are not contained in the allowable region, but that $v$ is, which is a contradiction.  Thus the $n$th step must be a Case St, Case A, Case A', Case A'', or Case D, so the next step can be Case Fi (it is worth noting here that this also leads to the converse of the Golden Lemma).  This completes the proof $\blacksquare$

	Following the above analysis, one can easily derive a generalized version of Theorem \ref{graph and snakes}:  The elements of the fundamental group of the closed surface $S$ with wordlength $n$ (according to any presentation described in the Main Theorem) are in bijective correspondence to directed paths from St to Fi in the graph $G_S$.

\subsection{The Markov Chain} \label{The Markov Chain}

Here we will describe in detail the Markov states in our Markov chain for the surface of genus 2.  We will use explicitly the presentation of the group $G=<a,b,c,d | abABcdCD=1>$.  Let $\Omega$ be the set of generators of $G$ along with their inverses.

Each state in the Markov chain represents a step in a bi-infinite snake.  This means that it represents the combination of the information of a region, and the directed edge or edges of that step.

Case A: This includes 8 different Markov states A$_k, k \in \Omega = {a,b,c,d,A,B,C,D}$.  Here, $k$ is simply the label of the edge traversed in this step.

Case B: This includes 8 different Markov states B$_{jk}$ where $j$ and $k \in \Omega$ are the edges traversed in this step with $j$ on the left relative to the direction of traversal.

Case C: This includes 8 different Markov states C$_{jk}$ where $j$ and $k \in \Omega$ are the edges traversed in this step with $j$ on the left relative to the direction of traversal.

Case C': This includes 8 different Markov states C'$_{jk}$ where $j$ and $k \in \Omega$ are the edges traversed in this step with $j$ on the left relative to the direction of traversal.

Case D: This includes 8 different Markov states D$_{jk}$ where $j$ and $k \in \Omega$ are the edges traversed in this step with $j$ on the left relative to the direction of traversal.

Case E: This includes 16 different Markov states E$^L_{ijk}$ and E$^R_{ijk}$ where $i$, $j$, and $k \in \Omega$ are the edges traversed in this step with $i$ on the left of $j$ which is to the left of $k$ relative to the direction of traversal.  The $L$ superscript is used when the $j$ edge starts at the start of the $i$ edge and ends at the end of the $k$ edge, and the $R$ superscript is used in the other case.

Case A': This includes 16 different Markov states A'$^L_k$ and A'$^R_k, k \in \Omega$. Here, $k$ is the edge traversed in this step.  The $L$ superscript is used when a wedge has been cut from the left side (relative to the direction of the edge) of the Case A region to make this Case A' region and the $R$ superscript is used in the other case.

Case A'': This includes 16 different Markov states A''$^L_k$ and A''$^R_k, k \in \Omega$. Here, $k$ is the edge traversed in this step.  The $L$ superscript is used when a wedge has been cut from the left side (relative to the direction of the edge) of the Case A' region to make this Case A'' region and the $R$ superscript is used in the other case.

Let $\mathbb{A}$, this set of 88 Markov states, be the alphabet for our Markov chain.  Let $\Pi$ be a map from $\mathbb{A}$ to the vertices of $G_2$ that simply removes the subscripts and superscripts.  The transition matrix for our Markov chain will now be described.  Choose $J, K \in \mathbb{A}$.  Then $J$ can be followed by $K$ if and only if there is a directed edge in $G_2$ from $\Pi(J)$ to $\Pi(K)$, and a choice of images in $\Gamma$ of $J$ and $K$ can be made (edges and corresponding regions) such that the edge(s) of $K$ start where the edge(s) of $J$ end and the region of $K$ is contained in that of $J$ as one of the specific subregions corresponding to the choice of an edge in $G_2$.
 
\subsection{Infinite Snakes and Their Uniqueness Properties}

\begin{defn}[Infinite Snake]\label{Infinite Snake}
An Infinite Snake is the image in the universal cover of an element of the Markov chain.
\end{defn}

\begin{lem}\label{4 consec}
An Infinite Snake never contains four or more consecutive edges of an octagon in $\Gamma$ without containing all 8.
\end{lem}

Proof: Assume there are are between 4 and 7 consecutive edges of a given octagon in $\Gamma$ but not an 8th.  The first 4 edges can be ordered 1 through 4 depending on their ordering of traversal.  The octagon cannot be left and then returned to without violating the Golden Lemma.  The second edge must be either a Case A' or a Case B (if it were a Case A then the previous edge couldn't be in the octagon, if it were a Case A'' then it must be at least the 3rd edge, and if it were a Case C, C', D, or E, then it would violate our assumption about the octagon being incomplete).  Case B is ruled out now becuase then the 3rd edge cannot be in the same octagon.  Thus the third edge must be a Case A''.  The corresponding region doesn't include the 4th edge so we are done $\blacksquare$

\begin{thm}\label{inverse}
An Infinite Snake has an inverse which is also an Infinite Snake.
\end{thm}

Proof: We will construct the inverse of the Infinite Snake and show that it is indeed the image of some element of the Markov chain.  Let $\mathbb{J}=...J_{-2}J_{-1}J_0J_1J_2...$ be our Markov chain.  We will make its inverse $\mathbb{K}=...K_{-2}K_{-1}K_0K_1K_2...$.  Let $e=\{e_i\}^\infty_{-\infty}$ and $R=\{R_i\}^\infty_{-\infty}$ be edges and regions of an Infinite Snake image of $\mathbb{J}$ where $e_i$ is(are) the edge(s) corresponding to the $i^{th}$ step in the Infinite Snake (the image of $J_i$), and $R_i$ is the corresponding region.

Define $f=\{f_i\}^\infty_{-\infty}$ such that $f_i = (e_{-i})^{-1}$ where inverse here simply means traversing the edge or edges in the opposite direction.  Whenever $\Pi(J_i)$ does not equal A, A', or A'', it is straightforward to deduce what $K_{-i}$ need be, as there is only one state that has the required set of edges to project to $f_{-i}$.  The following rules hold:

If $J_i =$ B$_{jk}$, then $K_{-i} =$ D$_{k^{-1}j^{-1}}$

If $J_i =$ D$_{jk}$, then $K_{-i} =$ B$_{k^{-1}j^{-1}}$

If $J_i =$ C$_{jk}$, then $K_{-i} =$ C'$_{k^{-1}j^{-1}}$

If $J_i =$ C'$_{jk}$, then $K_{-i} =$ C$_{k^{-1}j^{-1}}$

If $J_i =$ E$^{R(L)}_{jkl}$, then $K_{-i} =$ E$^{R(L)}_{l^{-1}k^{-1}j^{-1}}$

Now suppose $\Pi(J_i)$ equals A, A', or A''.  If no octagon of $\Gamma$ has both edge $e_i$ and any edge of $e_{i+1}$ as boundary componenets, then $K_{-i} =$ A$_{j^{-1}}$ where $e_i = j$.  If an octagon of $\Gamma$ has both edge $e_i$ and any edge of $e_{i+1}$ as boundary componenets, but no edge of $e_{i+2}$, then $K_{-i} =$ A'$^{R(L)}_{j^{-1}}$ where $e_i = j$ and the octagon is on the left(right) of the directed edge $e_i$.  Lastly, if an octagon of $\Gamma$ has edge $e_i$ and any edge of $e_{i+2}$ as boundary components, then $K_{-i} =$ A''$^{R(L)}_{j^{-1}}$ where $e_i = j$ and the octagon is on the left(right) of the directed edge $e_i$.

We must now verify that $\mathbb{K}$ as constructed satisfies the Markov chain transition conditions.  In particular, we will check that $K_{-i}$ can follow $K_{-i-1}$ for arbitrary $i$.  Clearly, as the edge or edges $e_{i+1}$ continue those of $e_i$, so must the edges $f_{-i}$ follow those of $f_{-i-1}$.  There are only finitely many possiblities for $J_i$ through $J_{i+4}$.  By the above lemma, this sufficiently determines what $K_{-i-1}$ and $K_{-i}$ must be.  In each of these cases, one can see there is a directed edge in $G_2$ from $\Pi(K_{-i-1})$ to $\Pi(K_{-i})$ wth the region corresponding to $f_{-i}$ being the appropriate subregion of $f_{-i-1} \blacksquare$ 

\begin{thm} \label{sufficient edges}
If $\mathbb{J}=...J_{-2}J_{-1}J_0J_1J_2...$ and $\mathbb{K}=...K_{-2}K_{-1}K_0K_1K_2...$ are two elements of the Markov chain that have Infinite Snake images with the same edges $e=\{e_i\}^\infty_{-\infty}$ where $e_i$ corresponds to both steps $J_i$ and $K_i$, then $\mathbb{J}$ and $\mathbb{K}$ are identical.
\end{thm}

Proof:  When $\Pi(J_i) \in \{$B, C, C', D, or E $\}$, $J_i = K_i$ trivially.  If not, one can use the above lemma to see that $\#\{j<i | $ an edge of $e_j$ is on the boundary of one of the same octagons as $e_i\}<3$.  From the definition of the transition matrix, one can deduce that this number must also be the number of "primes" following the A of $J_i$ and $J_k$.  The superscript of these states, if if present, is determined by the which side of $e_i$ the common octagon is on, and the subscript is determined by the label of $e_i\blacksquare$

The following theorem describe the sense in which Infinite Snakes are geometric limits of snakes.

\begin{thm} \label{geometric limit}
Any finite segment of an Infinite Snake can be extended (by a bounded amount) to a snake between two points in $\Gamma$.
\end{thm}

Proof:  Let $\mathbb{J}=...J_{-2}J_{-1}J_0J_1J_2...$ be an element of the Markov chain who has our Infinite Snake as an image.   Let $e=\{e_i\}^\infty_{-\infty}$  and $R=\{R_i\}^\infty_{-\infty}$ be edges and regions of the Infinite Snake image of $\mathbb{J}$ where $e_i$ is(are) the edge(s) corresponding to the $i^{th}$ step in the Infinite Snake (the image of $J_i$), and $R_i$ is the corresponding region.  Let $\{e_i|a\le i \le b\}$ be the segment of interest.

We must now find appropriate beginnings and ends of our snake.  If $\Pi(J_a)$ is B or A, then $e_a$ has one initial vertex which we can let be the start of the snake.  If $\Pi(J_a)$ is A' or A'', we can prepend the snake with a Case A or a Case A and Case A' respectively and then the initial vertex of the Case A can be the start of the snake.  If $\Pi(J_a)$ is C, C', D, or E, we can prepend as necessary until we prepend a Case C with a Case B, whose initial vertex can be the initial vertex of the snake.

Now to find an appropriate end vertex of our snake, if $\Pi(J_b)$ is D, A, A', or A'', then $e_b$ has one final vertex which we can let be the end of our snake.  If $\Pi(J_b)$ is B, C, C', or E, we can postpend the snake with steps until we reach a Case C' which can then be followed by a Case D.  Its final vertex will be the end of our snake.

Now that we have a start vertex and an end vertex, we can verify that what we have is indeed a snake between them.  We will do this by showing that we have a sequence of steps corresponding to a directed edgepath in $G_2$ from St to Fi.  Start at St.  The next vertex in $G_2$ will be A or B by construction.  The edge from St to A or B is chosen appropriately according to the lables of the initial edge or edges in the snake.  We continue making the directed edgepath by appropriately connecting the prepended Cases until we get to $\Pi(J_a)$.  Then, for $a \le i < b$ we can choose the edge in $G_2$ from $\Pi(J_i)$ to $\Pi(J_{i+1})$ because of the following: by definition of the transition matrix for our Markov chain, $R_{i+1}$ is one of the $m$ Case $\Pi(J_i+1)$ regions in the Case $\Pi(J_i)$ region, $R(i)$.  Each of these $m$ regions correspond precisely to one of the $m$ edges in $G_2$ from $\Pi(J_i)$ to $\Pi(J_{i+1})$ so we choose the appropriate one.  Then if necessary, we choose the appropriate edges to any postpended Cases and then from the final Case A, A', A'', or D to Fi.  We now have that our segment of the Infinite snake is a piece of the snake corresponding (up to deck tranformation of the initial vertex) to a directed path in $G_2$ from St to Fi $\blacksquare$

\begin{cor} \label{direction}
Given a geodesic $L$ of $N$ and an Infinite Snake, there is at most one direction in which edges of that Infinite Snake can cross $L$.
\end{cor}

Proof:  Let $\mathbb{J}=...J_{-2}J_{-1}J_0J_1J_2...$ be a Markov chain who has our Infinite Snake as an image.   Let $e=\{e_i\}^\infty_{-\infty}$ be edges of the Infinite Snake image of $\mathbb{J}$ where $e_i$ is(are) the edge(s) corresponding to the $i^{th}$ step in the snake (the image of $J_i$).  Suppose that, for some Geodesic $L$ of $N$, an edge of $e_a$ crosses $L$ in one direction, but an edge of $e_b$ crosses it in the other direction.  WLOG assume $a \le b$.  Then by the previous theorem, there is a snake between some vertices $v_0$ and $v_n$ containing the segment $\{e_i|a\le i \le b\}$ of the Infinite Snake.  $v_n$ has to be on one particular side of $L$.  Thus, one of the two traversing edges must be in violation of the Golden Lemma $\blacksquare$

\begin{thm} \label{limit points}
An Infinite Snake has exactly two limit points on the boundary of $\mathbb{H}^2$.
\end{thm}

Before we prove this, we will need the following lemma.

\begin{lem} \label{separating}
Given any two distinct points on the boundary of $\mathbb{H}^2$ there exists a geodesic $L$ of $N$ that separates the two.  The same holds for a point on the bondary and a point in the interior.
\end{lem}

Proof:  Consider the geodesic $l$ between the two points. $l$ has infinite hyperbolic length so it cannot be contained in the interior of any of the finite-sized hyperbolic octagons bounded by geodesics of $N$.  Therefore, it must cross one of these geodesics.  $l$ must then intersect some geodesic $L$ of $N$ transversely (if it intersected some geodesic of $N$ non-transversely, thin it would be a geodesic of $N$, each of which intersect infinitely many other geodesics of $N$ transversly).  Since $l$ and $L$ are geodesics that intersect transversly in hyperbolic space, they must intersect in exactly one point.  Therefore the two limit points of $l$ lie on opposite sides of $L$.  The proof of the second statement follows similarly $\blacksquare$

Proof of Theorem \ref{limit points}:   Let $\mathbb{J}=...J_{-2}J_{-1}J_0J_1J_2...$ be an element of the Markov chain who has our Infinite Snake as an image.   Let $e=\{e_i\}^\infty_{-\infty}$ be edges of the Infinite Snake image of $\mathbb{J}$ where $e_i$ is(are) the edge(s) corresponding to the $i^{th}$ step in the snake (the image of $J_i$).  For $0 \le i < \infty$ let $v_i$ be the the (left-most relative to the direction of the edges if there are multiple) final vertex of $e_i$.  The sequence $e=\{e_i\}^\infty_0$ must have accumulation points in the closure of the Poincare disk.  By the above theorem, since no snake ends two distinct steps on a given vertex (this would violate the shortest edgepath part of the definition), the vertices of $e$ must be unique.  Thus, since the only accumulation points of vertices of $\Gamma$ are in the boundary of the disk, the same must hold for $e$.

Now suppose that $b_1$ and $b_2$ are distinct accumulation points of $e$.  Then by the Lemma, there is a geodesic $L$ of $N$ separating them.  Since points in the sequence $e$ get arbitrarily close (in the Euclidian metric) to both $b_1$ and $b_2$, there must be $i < j < k$ such that $v_i$ and $v_k$ are on one side of $L$ but $v_j$ is on the other.  This leads to a violation of the above Corollary.  Thus there is a unique accumulation point and it must be a limit point.  We will call this the "forward limit point."

By constructing the inverse of the Infinite Snake via the previous theorem, one can find another limit point at the other end of the Infinite Snake.  We will call this the "backward limit point."  Suppose that they are the same point $b_0$.  Now via the above lemma construct $L'$ separating $v_0$ and $b_0$.  By a similar argument to the above paragraph, one can again construct a contradicition to the above Corollary.  This completes the proof $\blacksquare$

\begin{defn}[Proper Infinite Snake] \label{PIS}
A Proper Infinite Snake is an Infinite Snake such that all of its regions contain one of its limit points and the regions of its inverse all contain the other limit point.
\end{defn}

\begin{thm} \label{uniqueness}
Given any two distinct points $f$ and $b$ on the boundary of $\mathbb{H}^2$ that do not bound a geodesic of $N$, if $S$ and $S'$ are Proper Infinite Snakes with $f$ as forward and $b$ as backward limit points, then $S$ and $S'$ are the same. 
\end{thm}

Proof:  Let $\mathbb{J}=...J_{-2}J_{-1}J_0J_1J_2...$ be an element of our Markov chain with $S$ as image.   Let $e=\{e_i\}^\infty_{-\infty}$ be edges of the Infinite Snake image of $\mathbb{J}$ where $e_i$ is(are) the edge(s) corresponding to the $i^{th}$ step in the snake (the image of $J_i$).  We will divide this problem into three cases.

Case 1: There exists $i$ such that $\Pi(J_i)$ is B or D.  By considering inverses, we can assume it is D WLOG.  Let $R_i$ be the corresponding region for this step.  $f$ is clearly in $R_i$.  Consider the two geodesics that contain boundary rays of $R_i$.  Let $U$ be the open region between them.  Since $U$ contains the Case B region (call it $T$) corresponding to the inverse of step $i$, $b$ must also be in $U$.   By the Corollary, $S'$ must then be contained in $U$.  Let $v$ be the final vertex of $e_i$ and let $V$ be the octagon of $N$ containing it.  Since $V$ separates the parts of $U$ containing $f$ and $b$, $S'$ must pass through $V$ and thus contain $v$.  The only step with edges in $U$ that can have $v$ as a final vertex, $f$ contained in the corresponding region, and $b$ contained in the inverse step's region is the Case D step of $S$.  Thus $S'$ contains this step.  Since at every step, only one next step has a region that also contains $f$, the subsequent steps of $S'$ must be identical to those of $S$.  To determine what the preceding steps must be, one can use this same argument on the inverses of $S$ and $S'$ to first determine that the edges must be the same.  Using the process for constructing inverses of Infinite Snakes, we then uniquely construct the inverse of the inverse of the remaining part of $S'$.  The steps themselves (in particular whether an A, A', or A'' is present) are uniquely determined from the edges.

Case 2: There exists $i$ such that $\Pi(J_i)$ is E.  WLOG (by mirror symmetry) we can assume that the superscript of $J_i$ is $L$.  Let $R_i$ be the region for this step, and let $T$ be the inverse region for this step.  Consider the region $V$ made up of the two octagons of $N$ that each intersect two edges of $e_i$.  Let $L$ be the geodesic of $N$ on the common boundaries.  Let $U$ be the region bounded between the two geodesics of $N$ which contain the arc boundaries of these two octagons that are farthest from each other.  Since $f \in R_i \subset Int(U)$ and $b \in T \subset Int(U)$, $S'$ must also be in $U$ by the corollary.  As in Case 1, since $V$ separates the portions of $U$ containing $b$ and $f$, $S'$ must pass through $V$ and thus must contain at least one of the vertices of $\Gamma$ contained in $V$.  By rotational symmetry and by considering inverses, we can assume WLOG that $S'$ contains the vertex $v_r$ in the center of the right (relative to direction of $S$) octagon of $V$.  T  here are only two possible steps with edges in $U$ that can have $v_r$ as an initial vertex, $f$ contained in the corresponding region, and $b$ contained in an inverse step's region.  These are the Case E step of $S$ and a Case C step.  However based on the regions corresponding to the Case C step and its Case C' inverse, the only place $f$ and $b$ could be (and still be consistent with $S$) are on opposite ends of $L$, a contradiction.  Thus $S'$ contains the Case E step of $S$.  As in Case 1, we can then deduce further that $S$ is identical to $S'$.

Case 3: For all $i$, $\Pi(J_i) \in \{$A, A', A'' $\}$.  This can be divided into two sub-cases.

Case 3a: There exists $j$ such that $\Pi(J_j)$ is A'.  WLOG assume the supercript of $J_j$ is $L$, and let $\gamma \in N$ be the geodesic containing containing the left ray boundary of the corresponding region.  Let $v$ be the initial vertex of $e_j$ and let $V$ be the octagon of $N$ containing it.  $e_{j-1}$ must be the edge that enters $V$ immediately counterclockwise from where $e_j$ leaves it.  Now $e_{j-2}$ cannot be the edge crossing $\gamma$ and ending at the initial vertex of $e_{j-1}$, otherwise $\Pi(J_j)$ would be A''.  Thus the inverse edge of $e_{j-1}$ must have its region on the $V$ side of $\gamma$.   $b$ is in this region, and $f$, in the region for $e_i$, must be on the same side of $\gamma$.  Let $U$ be the region bounded by $\gamma$ and the six geodisics that contain boundary arcs of $V$ not crossed by $e_j$ or $e_{j-1}$.  $f$ and $b$ are in $U$, but are on opposite sides of $V$.  Thus, by the corollary, $S'$ must pass through $V$ and contain $v$.  The edges of $S'$ can be uniquely constructed outwardly from $v$ since at every vertex along the way, $f$ and $b$ are contained in the eminating regions of $S$ and its inverse and there is only one edge leading from the vertex into each of those.  The edges then determine the steps of $S'$ uniquely.

Case 3b:  There exists $k$ such that $\Pi(J_k)$ is A.  Let $v$ be the initial vertex of $e_k$ and let $V$ be the octagon of $N$ containing it.  Notice that $e_{k-1}$ and $e_k$ do not interesect adjacent edges of $V$.  Thus, the region $U$ bounded by the 6 geodesics containing the other edges of $V$ is disconnected when $V$ is removed.  The two disconnected regions must each contain either $b$ or $f$ since they contain the regions corresponding to $e_k$ and the inverse of $e_{k-1}$.  The argument proceeds as in Case 3a and we are done $\blacksquare$

\subsection{The Existence of Infinite Snakes via Algorithm}

We are trying to establish a correspondence between periodic elements of our Markov chain and conjugacy classes of $G$.  The previous result gives us a form of uniqueness of such an element.  Now we need to show the existence of such an element given a conjugacy class.  To this end, we will describe an algorithm, nearly identical to Lustig's \cite{lustig} but adapted to our language, for creating an Infinite Snake (that is an image of a Markov Chain) which has the same limit points as a lift of a given conjugacy class of $G$.  This will require an easier way to recognize Infinite Snakes.

First, we will describe topological or graphical structure of an Infinite Snake called a Topologic Snake.

\begin{defn}[Topological Snake] \label {TS}
A Topological Snake is a (finite, infinite, or bi-infinite) sequence of arcs and/or bubble chains.  A bubble chain is a sequence (finite, infinite, or bi-infinite) of bubbles (topologic circles) where each bubble shares a single proper topological interval with each neighboring bubble in the sequence, but no points with any other bubbles.  In each bubble chain, we assign a marked point to the part of the first (last) bubble disjoint from "boundary intervals", and we call this the entry (exit) point.  In case the bubble chain contains only one bubble, we require its entry point and exit point be distinct.  The arcs are directed, meaning they also have an entry point and an exit point (the endpoints).  The entry point of each bubble chain or arc is identified with the exit point of the preceeding bubble chain or arc.
\end{defn}

Now we will describe sufficient conditions for an embedding of a Topological Snake in $\Gamma$ to be an Infinite Snake.  We will call such an embedding a Geomtrical Snake.

\begin{defn}[Geometrical Snake] \label {GS}
A Geometrical Snake is an embedding in $\Gamma$ of a Topologal Snake that satisfies the following 4 Properties:

1.	Each bubble is mapped to an octagon of $\Gamma$.

2.	Each bubble must have two "special points" chosen from entry points, exit points, or endpoints of boundary intervals.  These two points must embed to opposite vertices of an octagon.  If the bubble has an entry point and/or exit point (if it is the first and/or last in its bubble chain), such point(s) must be the special point(s).

3.	Two bubbles in the same bubble chain cannot share a special point.

4.	If the embedding contains four or more edges of an octagon in $\Gamma$, it contains all eight.
\end{defn}

\begin{thm} \label {GS=IS}
A bi-infinite Geometrical Snake is an image of some element of the Markov chain and must also be an Infinite Snake.
\end{thm}

Proof: We will show how to construct a Markov chain whose image is a given Geometrical Snake, $S$.  First we partition the edges of $S$ into steps.  Each edge in the image of one of the arc sections of the Topological Snake will be in its own step of type A, A', or A''.  For edges in the images of bubbles, we consider the geodesic segment bewteen the special points of that bubble.  We group two edges of that bubble in the same step if they are mirror images of each other about that segment.  Properties 1 and 2 allow this to be clearly defined.  If no edge in one of these steps is a boundary interval, we have a unique choice of a B, C, C', or D step.  If an edge is a boundary interval, then it is grouped with its reflection in each of the the two adjacent octagons.  Property 3 ensures that these three edges form an E step.  As in the construction of inverse snakes, we can use Property 4 to distinguish between steps of type A, A', and A'', and then check that steps of types D, B, A, A', and A'' follow each other in ways consistent with the trasition matrix.  Steps of type B, C, C', D, and E follow each other in ways consistent with the transition matrix by construction and by Properties 1, 2, and 3 $\blacksquare$

It is now time to present an algorithm that converts a conjugacy class into a Geometrical Snake.  We start with a cyclic word in the generators of $G$.  The algorithm first uses the relation to shorten the word if possible, and then it adds other edges to include other loops in the conjugacy class of the same length.  The next theorem will prove that this algorithm always yeilds the properties characteristic of a Geometrical Snake.

Let $W=w_1w_2...w_l$ be our (nontrivial) cyclic word in the generators of $G$.  The indices here are understood modulo $l$, the current wordlength (which can shrink as we use relations).  The algorithm, wich will be proved to be finite, will consist of five steps.

Step 1.  Scan $W$ for instances where $w_i$ and $w_{i+1}$ are inverses and remove them.  Iterate this process until there are no more such instances.  Call the resulting cyclic word $W'$.

Step 2.  Scan $W'$ for any "long cycles" and replace them with the corresponding equivalent "short cycle".  A "long (short (half)) cycle" is a sequence of generators of $G$ that comprise more than (less than (exactly)) half of a permutation of the relation or its inverse.  Repeat Step 1 if necessary.  Then iterate this process until there are no more long cycles.  Call the resulting cyclic word $W''=w''_1w''_2...w''_{l''}$

Step 3.  Scan $W''$ for any half cycle $w''_{i}...w''_{i+3}$.  Replace it with the equivalent half cycle $v_1v_2v_3v_4$ if and only if $v_4w''_{i+4}w''_{i+5}w''_{i+6}$ is again a half cycle.  Repeat this procedure with the new half cycle and so on (the next theorem will prove that this process is finite).  Any time a long cycle or neighboring inverses appear, eliminate them as in Step 2.  Iterate this process until all half cycles have been checked and, when needed, replaced.  Call the resulting cyclic word $W'''$

Step 4.  Let $E=...e_{-1}e_0e_1...$ be an edgepath in $\Gamma$ that periodically reads off the cyclic word $W'''$.

Step 5.  Any time there are 4  consecutive sides of an octagon in $E$, add the other 4 and do this for the corresponding octagon in each of the other periods, creating a new Topological snake.  Iterate this process for each new Topological Snake until it cannot be done anymore.  Call the final collection of edges $GS$.

\begin{thm} \label {algorithm}
The above process yields a Geometrical Snake that has the same limit points as a lift of $W$.
\end{thm}

Before proving this theorem, we will need a lemma about the minimality of $E$.

\begin{lem} \label {locally minimal}
$E$ is locally minimal in the sense that any finite piece of it $e_ie_{i+1}...e_{j-1}e_j$ is a shortest edgepath.
\end{lem}

Proof: Steps 1 and 2 are clearly finite since every replacement involves a shortening.  In step 3, we need to insure that no replacement of a half cycle creates another replacement which creates another replacement and so on endlessly.  Note that when one replacent gives rise to another, both replacements must have the same orientation (both changing clockwise cycles to counter-clockwise or vise versa).  This must be the case since Steps 1 and 2 keep the word reduced.  After continuing such consecutive replacements once through the lenghth of the word, we come across  a cycle that has already been switched to the opposite orientation, so the procedure must halt.  When we change one half cycle, at most two new half cycles can be created, one before and one after.  If one is created before, it can be immediately confirmed that it does not need to be replaced.  If a half cycle is replaced because it creates a new one after, we can check that it does not need to be replaced again.  If the new succeeding half cycle does not need to be replaced, then we have a net decrease in the number that need to be checked.  Othewise we continue with the next one and the next one, but this can't go on forever, and when it stops, the net number that need to be checked decreases.  Thus, the net number of half cycles that need to be checked can never increase, but must eventually decrease, so the process terminates.

We will now introduce some terminology of Birman and Series \cite{birman-series} so that we can use a theorem of theirs to prove the lemma.  Two cycles $c_1$ and $c_2$ are called "consecutive" if they lie on the boundaries of different octagons in $\Gamma$ sharing a common edge $E$, such that $c_1$ ends on the endpoint of $E$ where $c_2$ begins, but neither cycle contains the other endpoint of $E$.  A sequence of consecutive cycles is called a "chain".  A "long chain" is a chain of $n>1$  cycles, $c_1,...,c_n$, where $c_1$ and $c_n$ are half cycles and all the rest have one less edge than these".  Such a chain of $3n+2$ can be replaced by the chain of length $3n$ running along the other side of the octagons, and this is exactly what Step 3 does.

A special case of Theorem 2.8(a) of Birman and Series \cite{birman-series} says that any edgepath that doesn't contain a long chain or long cycle is minimal.  Thus we only need to show that $E$ contains no long chains or long cycles to prove the lemma.  Suppose that $E$ does contain a long cycle.  Then $W'''$ must be at least as long as this cycle because the relation doesn't repeat symbols.  Thus the long cycle would be contained in a permutation of $W'''$, contradicting Step 2.  Now suppose $E$ contains a long chain.  Since it is finite, it must be contained in an edgpath reading off $(W''')^n$.  After cyclic permutation, we can assume WLOG that the chain starts at the beginning of $(W''')^n$.  The chain must be longer than $W'''$ or we'd have a contradiction of Step 3.  Thus we can infer that $W'''$ is one chain consisting of a half cycle followed cycles of one less length.  Step 3 would eliminate this scenario too $\blacksquare$

Now we prove Theorem \ref{algorithm}: $E = ...e_{-1}e_0e_1...$  must be an imbedded arc by the previous theorem so it is a Topological Snake.  It also vaccuously satisfies Properties 1, 2, and 3 of being a Geometrical Snake.  We must show that Step 5 is finite and that each addition of an octagon does not alter the facts that Properties 1, 2, and 3 are satisfied and that it is a Topological Snake.  Also, on the way to satisfying Property 4, we will ensure the weaker condition, which we call Property 4', that at every step, there are never 5-7 consecutive edges of an octagon without all 8 being present.  For ease of bookkeeping, when we complete octagons, we will group the new edges with their mirror images as in the proof of the previous theorem.  Let $F = ...f_{-1}f_0f_1...$ where $f_i$ contains $e_i$ and any new edge grouped with it.  It will follow inductively that for any new edge added to $f_i$, there will exist $m<\infty$ such that there is an edgepath of lenghth $m+1$ starting with the $e_{i-m}$ and ending with that new edge of $f_i$ and similarly there is an edgepath of length $m+1$ starting with this new edge and ending with $e{i+m}$.  Call this Property 5.  We will also maintain property called Property 6 that any time some $f_a$ contains two edges of an octagon, $F$ must contain all of them.  We will now break down the completion of octagons in Step 5 into three types based on the number of adjacent (sharing an edge) octagons already completed.

Type 0: we complete an octagon that shares no edges with any previously completed octagon.  To do this, we start with a half cycle $f_i...f_{i+3} = e_i...e_{i+3}$ that is the image of an arc or a piece of an arc of our Topological Snake.  The four new (all four are new by Property 4') edges add a bubble chain of one bubble to our Topological Snake.  They also are added simultaneously to ensure Property 6 remains satisfied.  This bubble is mapped to an octagon, consistent with Property 1.  The entry and exit point are opposite and can be chosen to be the special points consistent with Property 2.  Property 3 is still vacuously true for this bubble chain.  Property 5 will remain true trivially by simply using paths that go through the special points of this octagon.  Now, to see that we still have a Topological Snake, we must check that none of the three new vertices already lie on different edges of $F$.  If this were the case, we could create an edgepath using Property 5 that contradicts the minimality of $E$.  Thus the new bubble is properly embedded.  Similarly, Property 4' must still be satisfied.  Otherwise we could use Property 5 to yeild a contradiction of $E$ being minimal.

Type 1: we complete an octagon that shares exactly one edge with a previously completed octagon.   Start with a half cycle $d_j...d_{j+3}$ contained in $f_j...f_{j+3}$.  By Property 5 we must have that any half cycle must be contained in four consectuve pieces of $F$ to avoid contradicting the local minimality of $E$.  By Property 6 we have that it must span all four.  The edge bordering another completed octagon cannot be in $f_{j+1}$ or $f_{j+2}$ or we would have a contradiction of Property 2.  Thus $d_a=e_a=f_a$ for $a \in \{j+1,j+2\}$.  Thus, since by hypothesis only one edge of the halfcycle borders a previously completed octagon, we have either $d_a=e_a=f_a$ for $a=j$ or for $a=j+3$.  Assume the former is true (the other case is similar) so that $d_{j+3}$ will border the previously completed octagon and its initial vertex, being the entry point of the bubble chain, must be a special point of this octagon.  We now complete our current octagon.  Properties 6 and 1 immediately remains satisfied.  We have prepended the bubble chain containing the previously completed octagon with another bubble.  Let the endpoints of our half cycle being completed be the new octagon's special points.  This satisfies Property 2.  It also satisfies Property 3 since special points of the neighboring octagons are on opposite sides of $d_{j+3}$.  The other required Properties are checked to be satisfied as in Type 0 additions.

Type 2: we complete an octagon that shares exactly two edges with ah previously completed octagon.  As above, we start with a half cycle $d_k...d_{k+3}$ contained in $f_k...f_{k+3}$.  Again we must have $d_a=e_a=f_a$ for $a \in \{k+1,k+2\}$.  Thus $d_k$ and $d_{k+3}$ must be the edges bordering previously completed octagons.  This means that they have special points on the endpoints of $d_k$ and $d_{k+3}$ nearest each other.  Now we complete the octagon.  We eleminate an arc between two bubble chains, joining them as one bubble chain with the addition of a new bubble.  Let the endpoints of our half cycle being completed be the special points of the new octagon.  All of the required properties can be checked to be satisfied as in the Type 1 additions.

In each of the above completion Types, the cardinality of $\{i | 1 \le i \le l''', e_i = f_i \}$ strictly decreases.  Thus this process terminates.  Since it cannot continue, $F$ must finally satisfy Property 4.

Steps 1, 2, and 3 did not change the conjugacy class of $W$ so there is a lift of $W$ with the same limit points as $E$.  Since $GS$=$F$ contains $E$, it too has the same limit points$\blacksquare$

\subsection{Summary of Correspondence}

Most of the time, when we create an Infinite Snake from a conjugacy class, it is a Proper Infinite Snake, and therefore unique.  We will explore and classify the only instances where this is not the case.  This discussion will be analogous to Theorem 3.5 in Lustig's paper \cite{lustig}.

Let $W$ be a cyclic word in the generators of $G$.  Let $S$ be an Infinite Snake gotten from $W$ via the above algorithm.  Let $\mathbb{J}=...J_{-2}J_{-1}J_0J_1J_2...$ be an element of the Markov chain with $S$ as image.  Let $R_i$ be the region corresponding to $J_i$ for all $i$.  Suppose that $S$ is not a proper infinite snake.  WLOG, suppose that its forward limit point, $f$ is not contained in $R_k$ for some $k$.  Since $f$ is a limit point of $S$ the boundary of the Poincare disc, the edges of $S$ get arbitrarily close to it in the euclidian metric.  Now, if $f$ weren't in $\overline{R_i}$ for some $i$, then it also wouldn't be in the closures of all subsequent regions.  This cannot be the case since if $f$ were a nonzero distance from the regions in the euclidian metric, then the edges of $S$ closer to $f$ would grow arbitrarily far away from their regions in the hyperbolic metric, which contradicts the definition of our regions.  Thus, $f \in \overline{R_i}$ for all $i$.  The periodicity of $S$ and the nesting of the regions then dictate that $f$ is not in any of the regions, so it must be in their boundaries and they all must be open.  Type D regions cannot occur without Type C regions, so this implies that all regions are of Type A, A', or A''.  Let $\gamma \in N$ be the boundary of $R_k$ containing $f$ as a limit point.  Then $\gamma$ must bound all of the regions.  This describes an infinite chain of cycles of length 3.  Iterating backwards, one sees that $\gamma$'s other limit point is the other limit of $S$.  We will call such an Infinite Snake "exceptional".  Using this we can prove the following key theorem.

\begin{thm} \label {correspondence}
For $n>6$, there is an $n$-to-1 correspondence between elements of the Markov chain of period $n$ and primitive conjugacy classes of $G$ of wordlength $n$.
\end{thm}

Proof: With our group presentation, all primative conjugacy classes that yield exceptional Infinite Snakes via our algorithm have wordlength 6.  For primative conjugacy classes of wordlength greater than 6, the algorithm then yeilds a Proper Infinite Snake by the argument in the preceeding paragraph.  Choose a primative cyclic word $W$ of wordlength $n > 6$.  Since there is an exceptional Infinitie Snake with limit points corresponding to every geodesic of $N$, $W$'s Proper Infinite Snake, $S$ must not have both its limit points on the same geodesic of $N$.  Thus, up to deck transformation, it is the uinique Infinite Snake for the conjugacy class of $W$ by Theorem \ref{uniqueness}.  There are then exactly $n$ elements of the Markov chain of period $n$ corresponding to the different shifts that read off $S$.  Conversely, given an element of period $n$, one can read off the left-most edges of its Infinite Snake to obtain a corresponding cyclic word of length $n\blacksquare$

Let $\mathbb{J}_n$ be the set of all elements of our Markov chain that have period $n$.  Let $F_n$ be the set of conjugacy classes of $G$ of wordlength $n$.  Let $p_n$ be the pointwise map from $\mathbb{J}_n$ to $F_n$ described in the above proof.  Let $\lambda_n$ and $\mu_n$ be the uniform measures on $\mathbb{J}_n$ and $F_n$ respectively.  

\begin{thm} \label{pushforward}
The pushforward under $p_n$ of $\lambda_n$ approaches $\mu_n$ in total variation norm (TV) as $n$ approaches infinity.
\end{thm}

Proof: The number of elements of $\mathbb{J}_n$ grows exponentially with $n$, and a counting argument reveals that the number of these that aren't primative is on the order of the square root of them.  Thus the fraction of non primative elements of $\mathbb{J}_n$ approaches 0.  Since these then have negligible effect on measure $\lambda_n$, the previous theorem yields the desired result$\blacksquare$

\section{The Self-Intersection Number}

For $c \in F_n$ let $N(c)$ be the minimal number of self-intersections of a closed curve in the free homotopy class representing $c$.  This number can be computed as Lustig did it \cite{lustig}.

Using the map $p_n$, $N$ can also be defined on elements of $\mathbb{J}_n$.  As in the paper by Chas and Lalley \cite{chas-lalley}, we can express this value as a sum.  Choose $\textbf{X} \in \mathbb{J}_n$.  Let $\sigma^i\textbf{X}$ denote the $i$th cyclic permutation of $\textbf{X}$.  Construct via Lustig's algorithm \cite{lustig} $H:(\mathbb{J}_n\times \mathbb{J}_n)\rightarrow \{0,1\}$ that gives 1 if the inputs determine the first instance of a linked pair and 0 otherwise.  Then we can write
$$N(\textbf{X}) = \sum_{1\le i < j \le n} H(\sigma^i \textbf{X} ,\sigma^j\textbf{X}).$$ 
As in Chas and Lalley, $H$ can be written as $H(\textbf{X},\textbf{Y}) = \sum_{k=1}^n h_k(\textbf{X} ,\textbf{Y})$ where $h_k$ satisfies H0, H1, H2, and H3 of Chas and Lalley \cite{chas-lalley}, described below.

\subsection{The Spine of a Snake}

Before we describe in detail how to calculate $h_k$, we introduce the \emph{spine} of a snake which is equivalent to the "Normal Edgepaths" in used in Lustig's paper \cite{lustig}.  The \emph{spine} of a snake (or Infinite Snake, or Proper Infinite Snake) is a topological arc that runs along the center of the snake.  In arc-like regions of a snake, the spine is coincident with the arc.  When the spine reaches the entry point of a bubble chain, it travels along the geodesic from the entry point to the center of the octagon of the first bubble in the bubble chain.  If it is the only bubble in the bubble chain, the spine continues along the geodesic from the center of the octagon to the exit point.  When there are subsequent bubbles in the same chain, the spine travels along the geodesic from the center of the first octagon, and then to the next, until it reaches the center of the last bubble in the bubble chain.  Then it exits along the geodesic from the center of the last octagon to the exit point as in the case of a bubble chain with a single bubble.

We will construct a step-by-step map from a snake to its spine.  We will lose information about the Markov structure of the snake and the regions for each step, but we will keep exactly enough information to calculate self-intersection number via Lustig's algorithm \cite{lustig}.  The spine, as described above, is an embedded path in the graph $\Delta$ which we describe now.

\begin{defn}[$\Delta$] \label{Delta}
The graph $\Delta$, embedded in the hyperbolic plane, is the union of the graphs $\Gamma$ and $N$, together with geodesic edges between each vertex of $\Gamma$ and its 8 closest vertices of $N$.
\end{defn}

Now we construct the map $\phi$ from the alphabet $\mathbb{A}$ to a smaller alphabet $\Phi$ that translates graphically into a map from a snake in $\Gamma$ to its spine in $\Delta$: let $x,y,$ and$z$ be arbitrary elements in $\Omega$, the set of generators of $G$ along with their inverses.  Then we set

$\phi($A$_x)=\phi($A'$^L_x)=\phi($A'$^R_x)=\phi($A''$^L_x)=\phi($A''$^R_x)=\alpha_x$

$\phi($B$_{xy})=\beta_{xy}$,

$\phi($C$_{xy})=\phi($C'$_{xy})=\gamma$,

$\phi($D$_{xy})=\delta_{xy}$,

$\phi($E$^L_{xyz})=\epsilon_{y}$, and

$\phi($E$^R_{xyz})=\epsilon_{\bar{y}}$.

Here we will describe how $\phi$ translates graphically into a map from an Infinite Snake in $\Gamma$ to its spine in $\Delta$.
Suppose A$_x (x \in \Omega)$ is a step in a Markov chain and its Infinite Snake image is the edge $x_0 \in \Gamma$ with label $x$.  Then the spine of this step is the same edge in $\Delta$ and we label this edge (and any edge equivalent under deck transformation) $\alpha_x$, consistent with $\phi$.  The same holds for steps of the form A'$^L_x$, A'$^R_x$, A''$^L_x$, and A''$^R_x$.

Similarly, suppose B$_{xy} (x,y \in \Omega)$ is a step in an element of the Markov chain and its Infinite Snake image is the pair of edges $x_0$ and $y_0 \in \Gamma$ with labels $x$ and $y$.  Then the spine of this step is the edge in $\Delta$ from the vertex at the start of $x_0$ and $y_0$ to the vertex of $N$ at the center of the octagon of $\Gamma$ containing them.  We label this edge (and any edge equivalent under deck transformation) $\beta_{xy}$, consistent with $\phi$.

If C$_{xy} (x,y \in \Omega)$ is a step in an element of the Markov chain and its Infinite Snake image is the pair of edges $x_0$ and $y_0 \in \Gamma$ with labels $x$ and $y$, then the spine of this step is just the vertex in $N$ at the center of the octagon of $\Gamma$ containing the two edges.  We label this vertex (and any vertex equivalent under deck transformation) $\gamma$, consistent with $\phi$.  We do the same for any step of the form C'$_{xy}$.

Instead, if D$_{xy} (x,y \in \Omega)$ is a step in an element of the Markov chain and its Infinite Snake image is the pair of edges $x_0$ and $y_0 \in \Gamma$ with labels $x$ and $y$, then the spine of this step is the edge in $\Delta$ from the vertex of $N$ at the center of the octagon of $\Gamma$ containing $x_0$ and $y_0$ to the vertex of $\Gamma$ where they meet. We label this edge (and any edge equivalent under deck transformation) $\delta_{xy}$, consistent with $\phi$.

Now suppose E$^L_{xyz} (x, y, z \in \Omega)$ is a step in an element of the Markov chain and its Infinite Snake image is the set of edges $x_0$, $y_0$, and $z_0 \in \Gamma$ with labels $x$, $y$, and $z$.  Then the spine for this step is the edge in $N \subset \Delta$ from the vertex of $N$ at the center of the octagon of $\Gamma$ containing $y_0$ and $z_0$  to the vertex of $N$ at the center of the octagon of $\Gamma$ containing $x_0$ and $y_0$.  We label this edge (and any edge equivalent under deck transformation) $\epsilon_y$, consistent with $\phi$.

Finally suppose E$^R_{xyz} (x, y, z \in \Omega)$ is a step in an element of the Markov chain and its Infinite Snake image is the set of edges $x_0$, $y_0$, and $z_0 \in \Gamma$ with labels $x$, $y$, and $z$.  Then the spine for this step is the edge in $N \subset \Delta$ from the vertex of $N$ at the center of the octagon of $\Gamma$ containing $x_0$ and $y_0$  to the vertex of $N$ at the center of the octagon of $\Gamma$ containing $y_0$ and $z_0$.  We label this edge (and any edge equivalent under deck transformation) $\epsilon_{\bar{y}}$, consistent with $\phi$.

We can define inverses of the elements of $\Phi$ by pushing forward via $\phi$ the inverses described in Theorem \ref{inverse}.  This will give a well-defined map from $\Phi$ to itself, even though the inverse on single elements of $\mathbb{A}$ wasn't always well-defined.  In particula\r, we have $$\overline{\alpha_x}=\alpha_{\bar{x}}, \overline{\beta_{xy}}=\delta_{\bar{y} \bar{x}}, \overline{\gamma}=\gamma, \overline{\delta_{xy}}=\beta_{\bar{y} \bar{x}},\text{ and } \overline{\epsilon_y}=\epsilon_{\bar{y}}.$$
Notice that if an edge in $\Delta$ is labeled $\omega \in \Phi$, then the same edge, when traversed in the opposite direction, will be labeled $\overline{\omega}$.

\subsection{Calculating Linked Pairs} \label{Linked Pairs}

Now we can discuss the cyclic ordering of edges around the vertices of $\Delta$.  This will be used to identify "linked pairs" of Cohen and Lustig \cite{cohen-lustig} \cite{lustig}, corresponding to self-intersections.  Here and elsewhere we assume that $\Gamma$ is embedded the hyperbolic plane with the orientation that makes a counterclockwise traversal of the edges in an octagon read off the relation $abABcdCD$ of $G$ in order (up to cyclic permutation).  Then, the edges of $\Delta$ eminating from a vertex $v \in \Gamma \subset \Delta$ are, when read in counterclockwise order (up to cyclic permutation), $$O_1 = \alpha_{a} \beta_{da} \alpha_{d} \beta_{Cd} \alpha_C \beta_{DC} \alpha_D \beta_{cD} \alpha_c \beta_{bc} \alpha_b \beta_{Ab} \alpha_A \beta_{BA} \alpha_{B} \beta_{aB}.$$  Also, the edges of $\Delta$ eminating from a vertex $\gamma \in N \subset \Delta$ are, when read in counterclockwise order (up to cyclic permutation) $$O_2 = \delta_{AD} \epsilon_A \delta_{Ba} \epsilon_B \delta_{ab} \epsilon_{a} \delta_{bA} \epsilon_b \delta_{CB} \epsilon_C \delta_{Dc} \epsilon_D \delta_{cd} \epsilon_c \delta_{Dc} \epsilon_d.$$

With these orderings, we can create functions that assign an orientation to any cyclic word in the elements of $\Phi$.  Let $\xi$ be a cyclic word in $\Phi$.  For $i \in \{1,2\}$, set $o_i(\xi)$ equal to $1$ if the letters of $\xi$ are all in $O_i$ and occur in cyclic order in $O_i$, $-1$ if they are all in $O_i$ and occur in reverse cyclic order in $O_i$, and $0$ otherwise.

With these orientations, we can construct functions $u_k$ and $v_k$ to detect intersections on pieces of spines of length $k$.  Consider two (finite or infinite) strings $\boldsymbol{\zeta}=\zeta_1 \zeta_2...$ and $\boldsymbol{\xi}=\xi_1 \xi_2...$ of letters in $\Phi$.  We will now define $u_k$ on such a pair of strings for each $k \ge 1$.  First, unless the following conditions hold, we set $u_k(\boldsymbol{\zeta},\boldsymbol{\xi})=0$:\\
(a) each argument of $u_k$ has length at least $k$, and \\
(b) $\zeta_1 \neq \xi_1, \zeta_k \neq \xi_k$, and $\zeta_i=\xi_i$ for all $1<i<k$.

Given that (a) and (b) hold for a particular $\boldsymbol{\zeta}$, $\boldsymbol{\xi}$, and $k$, we define
$$
u_k(\boldsymbol{\zeta},\boldsymbol{\xi})=
\begin{cases}
1 & \text{if } k = 1 \text{, and } \zeta_1 \text{ and } \xi_1 \text{ are different elements of the set } \{\alpha_x, \epsilon_x\} \\
& \text{    for the same } x \in \Omega;\\
1 & \text{if } k = 2 \text{, and } o_1(\overline{\zeta_1} \text{ } \overline{\xi_1} \text { } \zeta_2 \text{ } \xi_2) = 1;\\
1 & \text{if } k = 3 \text{, and } o_1(\overline{\zeta_1} \text{ } \overline{\xi_1} \text { } \zeta_2) = o_1(\overline{\zeta_2} \text{ } \zeta_3 \text{ } \xi_3) \neq 0;\\
1 & \text{if } k = 4 \text{, and either } o_1(\overline{\zeta_1} \text{ } \overline{\xi_1} \text { } \zeta_2) = o_1(\overline{\zeta_3} \text{ } \zeta_4 \text{ } \xi_4) \neq 0 \text{, or} \\
& \text{    } o_2(\overline{\zeta_1} \text{ } \overline{\xi_1} \text { } \zeta_4 \text{ } \xi_4) = 1;\\
1 & \text{if } k \ge 5 \text{, and }  o_1(\overline{\zeta_1} \text{ } \overline{\xi_1} \text { } \zeta_2) + o_2(\overline{\zeta_1} \text{ } \overline{\xi_1} \text { } \zeta_4) = \\
& \text{    } o_1(\overline{\zeta_{k-1}} \text{ } \zeta_k \text { } \xi_k) + o_2(\overline{\zeta_{k-3}} \text{ } \zeta_k \text { } \xi_k) \ne 0; \text{ and}\\
0 & \text{otherwise.}
\end{cases}
$$

Now, for each $k$, we define $v_k(\boldsymbol{\zeta},\boldsymbol{\xi})=0$ unless each of its arguments has length at least $k$, in which case we set 
$$v_k(\boldsymbol{\zeta},\boldsymbol{\xi})=u_k(\zeta_1 \text{ } \zeta_2 \text{ ... } \zeta_k \text{, } \overline{\xi_k} \text{ } \overline{\xi_{k-1}} \text{ ... } \overline{\xi_1}).$$
$u_k$ and $v_k$ are defined for arbitrary finite or infinite strings of letters in $\Phi$.  We extend these definitions to \emph{doubly} infinite sequences $\boldsymbol{\zeta}=...\zeta_{-1} \zeta_0 \zeta_1...$, and $\boldsymbol{\xi}=...\xi_{-1} \xi_0 \xi_1...$ by adopting the convention that
$$u_k(\boldsymbol{\zeta},\boldsymbol{\xi})=u_k(\zeta_1 \text{ } \zeta_2 \text{ ... } \zeta_k \text{, } \xi_1 \text{ } \xi_2 \text{ ... } \xi_k) \text{, and}$$ 
$$v_k(\boldsymbol{\zeta},\boldsymbol{\xi})=v_k(\zeta_1 \text{ } \zeta_2 \text{ ... } \zeta_k \text{, } \xi_1 \text{ } \xi_2 \text{ ... } \xi_k).$$

Finally, we can define $h_k$.  If \textbf{X} and \textbf{Y} are (finite, infinite, or doubly infinite) strings of letters in $\mathbb{A}$ satisfying the Markov transition rules, we set $h_k(\textbf{X},\textbf{Y})=u_k(\phi(\textbf{X}),\phi(\textbf{Y})) + v_k(\phi(\textbf{X}),\phi(\textbf{Y}))$.  Here when we write $\phi$ acting on a string, it is to be interpreted as acting on all of the letters in the string.

\subsection{Properties of $h_k$ on $\lambda_n$}

We will need to show that $h_k$ satisfies 4 important properties of Chas and Lalley {chas-lalley}in order to make use of some of their results.  
They are:

(H0) Each function $h_k$ is symmetric.

(H1) There exists $C < \infty$ such that $|h_k| \le C$ for all $k \ge 1$.

(H2) For each $k \ge 1$ the function $h_k$ only depends on the first $k$ entries of its arguments.

(H3) There exist constants $C \le \infty$ and $0 < \theta < 1$ such that for all $n \ge k \ge 1$ and $0 \le i < j < n$,
$$E_{\lambda_n}|h_k(\sigma^i \mathbf{x},\sigma^j \mathbf{x})| \le C\theta^k.$$

The first three of these are clearly satisfied for our $h_k$ defined in Section \ref{Linked Pairs}.  Before we prove (H3), we will need some preliminary results.  We start with some properties of the transition matrix of our Markov chain.  Let $\mathbf{A} = \{a_{ij}| 1 \le i,j \le 88\}$ be the transition matrix of our topological Markov chain defined in Section \ref{The Markov Chain} on the alphabet $\mathbb{A}$ of 88 Markov states.  We have the following simple theorem.

\begin{thm} \label {irreducible}
The matrix $\mathbf{A}$ of $0$'s and $1$'s is irreducible and aperiodic.  Further, $\mathbf{A}^n$ is strictly positive for $n \ge 11$.
\end{thm}

Proof:  We will show that we can get from any state to A$_a$ and from A$_a$ to any state.  This will imply irreducibility.  From any state $J_1 \in \mathbb{A}$ such that $\Pi(J_1) \in \{$B, C, C', D, E$\}$ we can reach a state $J_2 \in \mathbb{A}$ such that $\Pi(J_2) =$D in 3 steps.  From any state $J_3 \in \mathbb{A}$ such that $\Pi(J_3) \in \{$A, A', A'', D$\}$ we can reach either A$_a$ or A$_c$ in one step and from A$_c$ we can reach A$_a$ in one step.  From A$_a$ we can get to A$_c$ in one step and from these we can reach any $J_4$ such that $\Pi(J_4)=$A.  From these we can reach any $J_5$ such that $\Pi(J_5) \in \{$A', B$\}$.  From these we can reach any $J_6$ such that $\Pi(J_6) \in \{$A'', C, C', D, E$\}$ in at most 3 steps.  This proves irreducibility.  The fact that we can get from A$_a$ to itself in one step then implies aperiodicity.  The specific statement for positivity follows from the fact that we can go from any state to A$_a$ in five or less steps, stay there for as many steps as we wish, and then go to any other state in 6 or less additional steps $\blacksquare$

Next we discuss a refinement of a general result from Section 8.5 of Walters \cite{walters} about the growth rate of periodic points of a Markov Chain.  Let $\tau: \mathbf{Y_M} \rightarrow \mathbf{Y_M}$ be a two sided topological Markov chain with irreducible, aperiodic $s$ by $s$ transition matrix $\mathbf{M}$.  Let $P_n$ denote the cardinality of the set $\mathbb{J}_n = \{\mathbf{y} \in \mathbf{Y_M} | \tau^n(\mathbf{y}) = \mathbf{y} \}$.  Let $\lambda_i (1 \le i \le s)$ be the eigenvalues of $\mathbf{M}$ and let $\Lambda$ be the eigenvalue that is positive and uniquely has the largest absolute value.  We have

\begin{thm} \label{growth}
$\operatorname{log}(P_n)/n$ approaches $\operatorname{log}\Lambda$ as $n$ approaches infinity.  This convergence is exponential in the following sense: there exists $C < \infty$ and $0 < \beta < 1$ such that $|P_n/\Lambda^n-1| \le C\beta^n$ for all $n$.
\end{thm}

Proof: A simple counting argument shows that $P_n$ equals the trace of $\mathbf{M}^n$ or $\sum_{i=1}^s \lambda_i^n$.  Let $\Lambda' < \Lambda$ be a number strictly greater in magnitude than all other eigenvalues of $\mathbf{M}$ (it exists by Theorem 1.1 of Seneta \cite{seneta}).  Combining this gives
$$ \left|\frac{\sum_{i=1}^s \lambda_i^n}{\Lambda^n}-1\right| \le (s-1) \Lambda'^n/ \Lambda^n.$$
Setting $C=s-1$ and $\beta = \Lambda'/\Lambda$ gives the desired exponential convergence and the result follows $\blacksquare$

We will also use the following fact about how well a spine determines a snake.

\begin{lem} \label{spine to snake}
A string of letters in $\Phi$ can have at most 16 preimages under $\phi$ that are portions of Infinite Snakes.
\end{lem}

Proof: Given the string $\zeta_1 \zeta_2...\zeta_n$, there are at most 16 preimages of $\zeta_1$ (this maximum occurs when $\zeta_1 = \gamma$.  Suppose that $J_1, J_2,$ and  $J_3$ are in $\mathbb{A}$ such that $J_2$ and $J_3$ can each immediately follow (according to the transition matrix $\mathbf{A}$) $J_1$, but $J_2 \neq J_3$.  Then it can be checked that $\phi(J_2) \neq \phi(J_3)$.  It inductively follows that a choice of preimage for $\zeta_1$ uniquely determines the rest of the preimage of the spine.  Thus there are at most 16 such preimages $\blacksquare$

The next lemma describes how a spine of a snake is "straight" and cannot turn back on itself or be it's own inverse.

\begin{lem} \label{straight}
Suppose $\zeta_1\zeta_2...\zeta_i$ be a portion of a spine such that $\zeta_j=\overline{\zeta_{i+1-j}}$ for all $j \le i$.  Then $i \le 2$.  
\end{lem}

Proof: If $i=3$, then $\zeta_2=\overline{\zeta_2}$, so $\zeta_2$ must be $\gamma$.  Exactly one of $\zeta_1$ or $\zeta_3$ must then also be $\gamma$, but this is impossible since they must be inverses of each other (this is checked by considering all possibilities for their preimages under $\phi$).  If $i=4$, then $\zeta_2$ and $\zeta_3$ must be inverses of each other.  The only way for this to be a portion of a spine is for them both to be $\gamma$.  If $x_1x_2x_3x_4$ is part of a snake such that $\phi(x_i) = \zeta_i$ where $\zeta_2 = \zeta_3 = \gamma$, then $\zeta_1$ and $\zeta_4$ cannot be inverses -- a contradiction.  If $i>5$ we would have one of the previous two cases (depending on whether $i$ were even or odd), so this too is impossible$\blacksquare$

Lastly, we will need a crude bound on the growth rate of strings of symbols in $\mathbb{A}$ satisfying the transition matrix $\mathbf{A}$.  Let $Q_n$ be the number of such strings of lenghth $n$.

\begin{lem} \label{strings}
$Q_n \leq 88^{11}P_n$.
\end{lem}

Proof: By Theorem \ref{irreducible}, any string of length $n-11$ can be made by removing the last $11$ letters from some element of $\mathbb{J}_n$.  Thus there are at most $P_n$ strings of length $n-11$.  The crude bound in the lemma comes from considering adding any $11$ elements of $\mathbb{A}$ to such a string $\blacksquare$  

We are now in position to prove that (H3) applies in our case.

\begin{thm} \label {H3}
There exist constants $C \le \infty$ and $0 < \theta < 1$ such that for all $n \ge k \ge 1$ and $0 \le i < j < n$,
$$E_{\lambda_n}|h_k(\sigma^i \mathbf{x},\sigma^j \mathbf{x})| \le C\theta^k.$$
\end{thm}

Proof: Since the measure $\lambda_n$ is invariant under cyclic shifts, and $h_k$ is symmetric in its inputs, it suffices to prove this in the case where $i=0$ and $j \le n/2$.  By property (H1), we simply need to show that the fraction of elements $\mathbf{x}$ of $\mathbb{J}_n$ such that $h_k(\mathbf{x},\tau^j \mathbf{x}) \neq 0$ falls off exponentially in $k$.  Let $P_n$ be the cardinality of $\mathbb{J}_n$.  Let $\Lambda$ be the maximal eigenvalue of $\mathbf{A}$.  By Theorems \ref{growth} and \ref{irreducible}, there exists $C_0 < \infty$ and $0 < \beta < 1$ such that for all $n$, 
\begin{equation} \label{growth approx}
|P_n/\Lambda^n-1| \le C_0\beta^n.
\end{equation}

Choose $n_0$ sufficiently large such that $C_0 \beta^{n_0} <1$.  For all $n > n_0$, rearragning (\ref{growth approx}) gives
\begin{equation} \label{P bound}
0<\Lambda^n(1-C_0\beta^{n_0})<\Lambda^n(1-C_0\beta^n)< P_n.
\end{equation}

We need to compute an upper bound on the number of elements $\mathbf{x}$ of $\mathbb{J}_n$ such that 
\begin{equation} \label{nonzero count}
h_k(\mathbf{x},\tau^j \mathbf{x}) \neq 0.  
\end{equation}

Suppose $\mathbf{x} = ...x_1x_2...x_n...$ satisfies (\ref{nonzero count}).  Then it is necessary that either
\begin{equation} \label{u neq}
\phi(x_2...x_{\lfloor(k-1)/2\rfloor})=\phi(x_{j+2}...x_{j+\lfloor(k-1)/2\rfloor}),
\end{equation}

or
\begin{equation} \label{v neq}
\phi(x_2...x_{\lfloor(k-1)/2\rfloor})=\overline{\phi(x_{j+k-1})} \text{ } \overline{\phi(x_{j+k-2})}...\overline{\phi(x_{j+k+1-\lfloor(k-1)/2\rfloor})}.
\end{equation}

Assume $\mathbf{x}$ satisfies (\ref{u neq}).  Then $\phi(x_2...x_{\lfloor(k-1)/2\rfloor})$ is determined uniquely by the $n-\lfloor(k-1)/2\rfloor$ symbols $\phi(x_{\lfloor(k-1)/2\rfloor+1})...\phi(x_n)$ which are determined by the string $x_{\lfloor(k-1)/2\rfloor+1}...x_n$ of length $n-\lfloor(k-1)/2\rfloor < n-k/2 + 1$.  By Lemma \ref{strings} and (\ref{growth approx}), there are at most $88^{11}\Lambda^{n-k/2+1}(1+C_0\beta^{n-k/2+1}) < C_1\Lambda^{n-k/2+1}$ such strings, where $C_1 = 88^{11}(C_0+1)$.  Each of these yield at most $16$ values for $x_2...x_{\lfloor(k-1)/2\rfloor}$ according to Lemma \ref{spine to snake}.  If we include the $88$ values $x_1$ can have, then there are maximum of $16(88)C_1\Lambda^{n-k/2+1}$ values of $\mathbf{x} \in \mathbb{J}_n$ satisfying (\ref{u neq}).

Now assume that $\mathbf{x}$ satisfies (\ref{v neq}).  By Lemma \ref{straight}, $j+k \le n+4$.  Then $\phi(x_5...x_{\lfloor(k-1)/2\rfloor})$ is determined uniquely by the $n-\lfloor(k-1)/2\rfloor$ symbols $\phi(x_{\lfloor(k-1)/2\rfloor+1})...\phi(x_n)$ which are determined by the string $x_{\lfloor(k-1)/2\rfloor+1}...x_n$ of length $n-\lfloor(k-1)/2\rfloor < n-k/2 + 1$.  As in the previous paragraph, there are at most $C_1\Lambda^{n-k/2+1}$ such strings.  Each yeilds at most $16$ values for $x_5...x_{\lfloor(k-1)/2\rfloor}$ according to Lemma \ref{spine to snake}.  If we include the $88$ values each of $x_1, x_2, x_3,$ and $x_4$ can take, then there are a maximium of $16(88^4)C_1\Lambda^{n-k/2+1}$ values of $\mathbf{x} \in \mathbb{J}_n$ satisfying (\ref{v neq}).

Combining the results of the previous two paragraphs, there are $16(88^4+88)C_1\Lambda^{n-k/2+1}$ values of $\mathbf{x} \in \mathbb{J}_n$ satisfying  (\ref{nonzero count}).  Setting $C_2 = 16(88^4+88)C_1\Lambda/(1-C_0\beta^{n_0})$ and using (\ref{P bound}) we have the upper bound

\begin{equation} \label {E lambda bound}
E_{\lambda_n}|h_k(\mathbf{x},\sigma^j \mathbf{x})| \le \frac{C_2\Lambda^{n-k/2}}{\Lambda^n}=C_2\Lambda^{-k/2}
\end{equation}

when $n \ge n_0$.  When $n < n_0$ we still have $E_{\lambda_n}|h_k(\mathbf{x},\sigma^j \mathbf{x})| \le 1$.  Thus, setting $C=\max(C_2,\Lambda^{n_0k/2})$ and $\theta=\Lambda^{-1/2}$ gives the desired result $\blacksquare$

\section{Statistics}

In this section, we state and prove the main theorem.

\begin{thm}[Main Theorem] \label{Main}
Let $N_n$ denote the random variable obtained by evaluating $N$ on a uniformly chosen $\alpha \in F_n$.  Then there exists $\kappa$ and $\sigma^2 \ge 0$ such that, as $n \rightarrow \infty$,
$$\frac{N_n-\kappa n^2}{n^{3/2}} \Longrightarrow Normal(0,\sigma)$$
in the weak topology.
\end{thm}

Theorem \ref{pushforward} implies that the distribution of $$f_n = \frac{\sum_{i=1}^n (\sum_{j=i+1}^n(\sum_k^\infty( h_k(\sigma^i \textbf{X} ,\sigma^j\textbf{X}))))-\kappa n^2}{n^{3/2}}$$ on $\mu_n$ and$\lambda_n$ approach each other.  Thus, to prove the Main Theorem, we only need to prove that the distribution of $f_n$ on $\lambda_n$ approaches a Gaussian.  To this end, we will prove (Theorem \ref{f and g}) that it is weakly equvialent to a normalized U-statistic, which according to Theorem 5.1 of Chas and Lalley \cite{chas-lalley} is Gaussian.  This will prove the Main Theorem.

\subsection{The Parry Measure}

Before we define the U-statistic, we need to define the Parry measure on our topological Markov chain.  By Theorem \ref{irreducible} and Theorem 1.1 of Seneta \cite{seneta}, $\mathbf{A}$ has a positive eigenvalue (of uniquely maximal absolute value) $\lambda$ with positive right and left eigenvectors.  Let $u=(u_1,...,u_{88})$ be a stricty positive left eigenvector and $v=(v_1,...,v_{88})$ be a strictly positive right eigenvector of $\mathbf{A}$ such that $\sum_{i = 1}^{88} u_iv_i=1$.  Then let $\pi_i=u_iv_i$ and $q_{ij}=a_{ij}v_j/(\lambda v_i)$.  Then, as in Walters \cite{walters}, the stochastic matrix $Q=\{q_{ij}\}$ and the stationary probablity vector $\pi=\{\pi_i\}$ give our Markov chain the Parry measure $\nu$.

The Parry measure has some important properties.  It is the unique measure with maximal entropy on our Markov chain $\mathbf{X_A}$ according to Theorem 8.10 of Walters \cite{walters}.  This implies according to Theorem 8.17 of Walters \cite{walters} that it "describes the distribution of the periodic points" of $\tau$.  This means that for a measurable set $S \subset \mathbf{X_A}$
\begin{equation} \label {describes}
\lim_{n \to \infty} \frac{\#\{\mathbf{x} \in \mathbb{J}_n|\mathbf{x} \in S\}}{P_n} = \lim_{n \to \infty} (\lambda_n(\mathbb{J}_n \bigcap S)) = \nu (S).
\end{equation}

We will use this result to produce an analogue of Theorem \ref{H3} for the Parry measure $\nu$.

\begin{thm} \label{H3 nu}
There exist constants $C \le \infty$ and $0 < \theta < 1$ such that for all $k \ge 1$ and $0 \le i < j$,
$$E_{\nu}|h_k(\tau^i \mathbf{x},\tau^j \mathbf{x})| \le C\theta^k.$$
\end{thm}

Proof: Fix $i, j,$ and $k$.  Let $C$ and $\theta$ be as in Theorem \ref{H3}.  Let $S = \{\mathbf{x} \in \mathbf{X} | h_k(\tau^i \mathbf{x},\tau^j \mathbf{x})=1\}$.  Note that $S$ is measureable since it is a finitely determined cylinder set.  Since the range of $h_k$ is $\{0,1\}$, 
\begin{equation} \label{E nu is nu}
E_{\nu}|h_k(\tau^i \mathbf{x},\tau^j \mathbf{x})|=\nu (S) \text{, and}
\end{equation}
\begin{equation} \label{E lambda is lambda}
E_{\lambda_n}|h_k(\sigma^i\mathbf{x},\sigma^j\mathbf{x})|=\lambda_n\{\mathbf{x} \in \mathbb{J}_n | h_k(\sigma^i\mathbf{x},\sigma^j\mathbf{x})=1\}
\end{equation}
For $\mathbf{x} \in \mathbb{J}_n \subset \mathbf{X_A}$, $\sigma \mathbf{x} = \tau \mathbf{x}$, so (\ref{E lambda is lambda}) becomes
\begin{equation} \label{E lambda is JS}
E_{\lambda_n}|h_k(\sigma^i\mathbf{x},\sigma^j\mathbf{x})|=\lambda_n(\mathbb{J}_n \bigcap S).
\end{equation}
By (\ref{E lambda is JS}) and Theorem \ref{H3}, we have bound (indepent of $n$)
\begin{equation} \label{JS bound}
\lambda_n(\mathbb{J}_n \bigcap S) \le C\theta^k.
\end{equation}
Finally, taking the limit on both sides of (\ref{JS bound}), and applying (\ref{describes}) explicitly gives $\nu(S) \le C\theta^k$, so (\ref{E nu is nu}) gives the desired result $\blacksquare$

We will also need a theorem relating a substring measure to the Parry measure.  Fix $0 < \delta < 1/4$ and set $m=m(n)=\lfloor n^\delta \rfloor$.  Let $\lambda_{n,m}$ denote the distribution of the random substring $x_1x_2...x_{n-m}$ of $x_1x_2...x_n \in \mathbb{J}_n.$  Also let $nu_{n,m}$ denote the distribution of the random substring $x_1x_2...x_{n-m}$ of $...x_1x_2...x_{n-m}... \in \mathbf{X_A}.$  Then

\begin{thm} \label {lambda n m}
For $n$ sufficiently large, $\lambda_{n,m}$ is absolutely continuous with respect to $\nu_{n,m}$, and the Radon-Nikodym derivative is sandwiched between bounds that tend to $1$ so that the two measures agree in total variation norm.
\end{thm}

Proof:  When $n \ge 11^{1/\delta}$, Theorem \ref{irreducible} guarantees absolute continuity because for any string $...x_1x_2...x_{n-m}... \in \mathbf{X_A}$, there exists at least one element $x_1x_2...x_n$ in $\mathbb{J}_n$ since $m$ will be at least $11$.  Now, at the substring $x_1x_2...x_{n-m}$,
$$\frac{d\lambda_{n,m}}{d\nu_{n,m}}=\frac{\sum\limits_{x_{n-m+1},x_{n-m+2},...,x_n \in \mathbb{A}} a_{x_1,x_2}a_{x_2,x_3}...a_{x_{n-m-1},x_{n-m}}a_{x_{n-m},x_{n-m+1}}...a_{x_n,x_1}/P_n}{u_{x_1}v_{x_1}\prod\limits_{i=1}^{n-m-1}(a_{x_i,x_{i+1}}v_{x_i+1}/(\Lambda u_{x_i}))}$$

$$=\frac{a_{x_1,x_2}a_{x_2,x_3}...a_{x_{n-m-1},x_{n-m}} (\mathbf{A}^{m+1})_{x_{n-m},x_1}/P_n}{a_{x_1,x_2}a_{x_2,x_3}...a_{x_{n-m-1},x_{n-m}}u_{x_1}v_{x_{n-m}}/\Lambda^{n-m-1}} = \frac{(\mathbf{A}^{m+1})_{x_{n-m},x_1}/P_n}{u_{x_1}v_{x_{n-m}}/\Lambda^{n-m-1}} $$
According to Theorem 1.2 of Seneta \cite{seneta}, $(\mathbf{A}^{r})_{x_i,x_j} \longrightarrow \Lambda^r u_{x_j}v_{x_i}$ as $r$ goes to infinity.  We can use this and the bounds from Theorem \ref{growth} to force the Radon-Nikodym derivitive bounds toward $1$ as follows.
$$\frac{d\lambda_{n,m}}{d\nu_{n,m}}=\frac{(\mathbf{A}^{m+1})_{x_{n-m},x_1}/P_n}{u_{x_1}v_{x_{n-m}}/\Lambda^{n-m-1}} \longrightarrow \frac{\Lambda^{m+1} u_{x_1}v_{x_{n-m}}/\Lambda^n}{u_{x_1}v_{x_{n-m}}/\Lambda^{n-m-1}} = 1$$
as $n \rightarrow \infty$ independently of the choice of $x_i$'s $\blacksquare$

\subsection{The U-Statistic}

We will describe an example of a second order U-statistic as in Chas and Lalley \cite{chas-lalley}.  When normalized, it will apprach a Gaussian by Theorem 5.1 of their paper and we will show it also weakly approaches $f_n$.  Call our Markov chain space $X_{\mathbf{A}}$.  Assign it the Parry measure as above and let $\tau$ be the shift map on the sequence space $\mathbb{A}^{\mathbb{N}}$.  For $\mathbf{X}$ chosen from this stationary, aperiodic, irreducible Markov chain with measure $\nu$, define the following second order U-statistic $$N'_n(\mathbf{X})=\sum_{i=1}^n(\sum_{j=i+1}^n(\sum_{k=1}^{\infty}(h_k(\tau^i\mathbf{X},\tau^j\mathbf{X})))).$$

Let $f_n$ be the distribution of $(N(c)-n^2\kappa)/n^{3/2}$ on $\mu_n$ (the uniform measure on the set of conjugacy classes $F_n$ of length n).  Let $g_n$ be the distribution of $(N'_n(\mathbf{x})-\kappa n^2)/n^{3/2}$ on the Parry measure $\nu$.  This is our "normalized" version of the U-statistic.  Before we relate this U-statistic to the intersection number, we will need an elementary result about distributions.

\begin{lem} \label {analysis}
Suppose $X_n$ and $Y_n$ are two infinite sequences of random variables that are all defined on a common probability space.  Let $c_n$ be an infinite sequence of scalars.  Choose $r>0$.  Let $F_n$ and $G_n$ denote the distributions of $X_n$ and $Y_n$, respectively.  Then,
$$\| F_n -G_n\|_{TV} \longrightarrow 0$$
implies that
$$(X_n-Y_n)/n^r \Longrightarrow 0,$$
which itself implies that
$$\frac{Y_n-a_n}{n^r} \Longrightarrow F \text{ if and only if } \frac{X_n-a_n}{n^r} \Longrightarrow F.$$
\end{lem}

The following theorem and proof closely follow Proposition 4.7 of Chas and Lalley \cite{chas-lalley}.

\begin{thm} \label{f and g}
For any metric $\varrho$ which induces the topology of weak convergence on probability measures,
\begin{equation} \label{varrho f and g}
\lim_{n \to \infty} \varrho(f_n, g_n) = 0, 
\end{equation}
so $f_n \Longrightarrow Normal(0,\sigma)$ if and only if $g_n \Longrightarrow Normal(0, \sigma)$.
\end{thm}

Proof: Let $f'_n$ be the distribution of $(N(\mathbf{x})-n^2\kappa)/n^{3/2}$ under the uniform measure $\lambda_n$ on $\mathbb{J}_n$.  By Theorem \ref{pushforward}, Lemma \ref{analysis}, and the fact that total variation norm is never increased by a mapping, $f_n$ and  $f'_n$ approach each other in total variation distance.  Thus, it suffices to prove (\ref{varrho f and g}) with $f'_n$ replacing $f_n$.  We will do this by truncating the formulas for the intersection number and U-statistic,

\begin{equation} \label{N}
N(\mathbf{x})=\sum_{i=1}^n(\sum_{j=i+1}^n(\sum_{k=1}^n(h_k(\sigma^i\mathbf{x},\sigma^j\mathbf{x})))) \text{ and}
\end{equation}
\begin{equation} \label{Nn}
N'_n(\mathbf{X})=\sum_{i=1}^n(\sum_{j=i+1}^n(\sum_{k=1}^{\infty}(h_k(\tau^i\mathbf{X},\tau^j\mathbf{X})))),
\end{equation}
twice each and then comparing the results using machienery similar to that in Walter's proof of (\ref{describes}).

Fix  $0 < \delta < 1/4$ and set $m=m(n)=\lfloor n^\delta \rfloor$.  By Theorems \ref{H3} and \ref{H3 nu}, 
$$E_{\lambda_n}( \sum_{i=1}^n(\sum_{j=i+1}^n(\sum_{k=m(n)}^n|h_k(\sigma^i\mathbf{x},\sigma^j\mathbf{x})|))) \le Cn^2\theta^{m(n)}, \text{ and}$$
$$E_{\nu}( \sum_{i=1}^n(\sum_{j=i+1}^n(\sum_{k=m(n)}^{\infty} |h_k(\tau^i\mathbf{x},\tau^j\mathbf{x})|))) \le Cn^2\theta^{m(n)}.$$
These upper bounds decrease to $0$ as $n\rightarrow \infty$.  Thus, by Markov's inequality, both of these sums converge weakly to $0$.  According to Lemma \ref{analysis}, this implies that
$$\lim_{n \to \infty} \varrho(f'_n, f'^A_n) = 0 \text{, and } \lim_{n \to \infty} \varrho(g_n, g^A_n) = 0.$$
Where $f'^A_n$ and $g^A_n$ are gotten from $f'_n$ and $g_n$ respectively by replacing the bounds on the inner sums of (\ref {N}) and (\ref{Nn}) by $1 \le k < m(n)$.  Therefore, it suffices to prove
$$\lim_{n \to \infty} \varrho(f'^A_n, g^A_n) = 0.$$

Notice that $|\{i,j,k| n-2m(n)< j < n, 1 \le i < j, 1 \le k < m(n)\}| \le 2(n)(m(n))^2 = 2n^{1+2\delta}$.  Since this grows much slower than $n^{3/2}$, we can eliminate this index set from the summations describing $f'^A_n$ and $g^A_n$ with negligible effect (each term is at most $1$).  Specifically, let $f'^B_n$ and $g^B_n$ be the distributions of the second truncations of $f'_n$ and $g'_n$ gotten by replacing the limits of summation in (\ref {N}) and (\ref{Nn}) by $1 \le i < j < n-2m(n)$ and $1 \le k < m(n)$.  Then by Lemma \ref{analysis} it suffices to prove 
$$\lim_{n \to \infty} \varrho(f'^B_n, g^B_n) = 0.$$
Here, when $i < j < n- 2m(n)$ and $k < m(n)$, $h_k(\sigma^i\mathbf{x},\sigma^j\mathbf{x})$ and $h_k(\tau^i\mathbf{X},\tau^j\mathbf{X})$ only depend on the first $n-m(n)$ entries of $\mathbf{x}$.  Thus $f'^B_n$ and $g^B_n$ are the same function on the two different measures $\lambda_{n,m}$ and $\nu_{n,m}$.  Finally, by Theorem \ref{lambda n m} and since total variation distance is never increased under a mapping, we have the desired result $\blacksquare$

\section{Positivity of Variance}

Here we prove that the $\sigma^2$ of Theorem \ref{Main} is strictly positive.  Moreover, the Gaussian is not degenerate.

\begin{thm} \label{variance}
$\sigma^2 > 0$ where $\sigma^2$ is the variance of the Gaussian in Theorem \ref{Main}.
\end{thm}

Before proving this, we will need to introduce the Hoeffding projection.

\subsection{Hoeffding Projection}
Given a kernel such as $h_k$ defined above, it's Hoeffding projection is defined  as $$H_k(\mathbf{z})=E_{\nu}h_k(\mathbf{z},\mathbf{Z}).$$
This is the expected value of $h_k$ when one of its inputs is fixed and the other is chosen randomly according to the parry measure.  We will now state a property of the decay rate of $H_k$.

\begin{thm} \label{Hoeffding decay}
$H_k(\mathbf{z})$ falls off exponentially with $k$ uniformly.  In other words, there exists constants $C \le \infty$ and $0 < \theta < 1$ such that for all $k$ and for all $\mathbf{z}$, $H_k(\mathbf{z}) \le C\theta^k$.
\end{thm}

Proof:  First we observe that there exists constants $C \le \infty$ and $0 < \theta < 1$ such that for all $k, n$ and $\mathbf{z}$, $E_{\lambda_n}(h_k(\mathbf{z},\mathbf{Z})) \le C\theta^k$.  This is checked routinely using the method of the proof of Theorem \ref{H3}.  The theorem then follows from an application of the proof of Theorem \ref{H3 nu} $\blacksquare$

Now, let $H=\sum_{k=1}^{\infty}H_k$, which must be finite by Lemma \ref{Hoeffding decay}.  We also have the following nice result about the continuity of $H$.

\begin{cor} \label{Holder}
H is Holder Continuous.  In other words, there exists $C_0 < \infty$ and $0<\theta_0<1$ such that for any $\mathbf{x} = ...x_{-1}x_0x_1...$ and $\mathbf{y} = ...y_{-1}y_0y_1... \in \mathbf{X_A}$, if $x_i$ = $y_i$ when $|i| \le n$ then $|H(\mathbf{x}) - H(\mathbf{y})| \le C_0\theta_0^n$.
\end{cor}

Proof: Suppose $\mathbf{x} = ...x_{-1}x_0x_1...$ and $\mathbf{y} = ...y_{-1}y_0y_1... \in \mathbf{X_A}$ where $x_i$ = $y_i$ when $|i| \le n$.  Then $H_k(\mathbf{x})=H_k(\mathbf{y})$ for $k \le n$ since $h_k$ and thus $H_k$ only depend on the first $k$ entries.  Thus,
$$|H(\mathbf{x})-H(\mathbf{y})| = \sum_{k=n+1}^{\infty}|H_k(\mathbf{x})-H_k(\mathbf{y})|$$
Since $H_k$ is nonnegative we can use Theorem \ref{Hoeffding decay} to bound the above expression by
$$\sum_{k=n+1}^{\infty}C\theta^k=\frac{C\theta}{1-\theta}\theta^n$$
Setting $C_0 = C\theta/(1-\theta)$ and $\theta_0=\theta$ gives the desired result $\blacksquare$

\subsection{Comparing Measures}

Before we proceed, we will nee a few tools for comparing parry measures of cylinder sets.  Call $s=J_1J_2...J_n$ a $string$ of length $n$ if the $J_i$'s are in $\mathbb{A}$ such that the $J_iJ_{i+1}$ satisfies the transition matrix $\mathbf{A}$.  Define the cylinder set of $s$ as $cyl(s)=\{...K_{-1}K_0K_1...\in X_{\mathbf{A}}|K_i=J_i$ for $1 \le i \le n\}$.

\begin{thm} \label{compare}
Suppose $j=J_1...J_n$ and $k=K_1...K_n$ are strings of length $n$ such that for all $i \in \mathbb{N}$, the number $J^{i+}$ of strings of length $i+n$ of the form  $jJ_{n+1}...J_{n+i}$ is the same as the number of strings $K^{i+}$ of the form $kK_{n+1}...K_{n+i}$, and the number of strings $J^{i-}$ of length $i+n$ of the form $J_{1-i}...J_0j$ is the same as the number of strings $K^{i-}$ of the form $K_{1-i}...K_0k$.  Then $\nu(cyl(j))=\nu(cyl(k))$.
\end{thm}

Proof: We will again use the result from Seneta \cite{seneta}, $(\mathbf{A}^{i})_{x,y} \longrightarrow \Lambda^i u_{y}v_{x}$ as $i$ goes to infinity.  In other words,
\begin{equation} \label{seneta eq}
\lim_{i \to \infty}\frac{(\mathbf{A}^{i})_{x,y}}{\Lambda^i u_{y}v_{x}}=1.
\end{equation}
Here $\mathbf{A}$ is our irreducible, aperiodic transition matrix, $x$ and $y$ are in $\mathbb{A}$, and $u$ and $v$ are the left and right eigenvectors of the leading eigenvalue $\Lambda$ of $\mathbf{A}$.  Now a counting argument gives

\begin{equation} \label{count strings}
\sum_{l \in \mathbb{A}}(\mathbf{A}^i)_{J_n,l} = J^{i+} = K^{i+} = \sum_{l \in \mathbb{A}}(\mathbf{A}^i)_{K_n,l}.
\end{equation}
By applying (\ref{seneta eq}) for all $l \in \mathbb{A}$ we can deduce
$$\lim_{i \to \infty}\frac{\sum_{l \in \mathbb{A}}((\mathbf{A}^{i})_{J_n,l})}{\sum_{l \in \mathbb{A}}(\Lambda^i u_{l}v_{J_n})}= 1 =\lim_{i \to \infty}\frac{\sum_{l \in \mathbb{A}}((\mathbf{A}^{i})_{K_n,l})}{\sum_{l \in \mathbb{A}}(\Lambda^i u_{l}v_{K_n})},$$
whence (\ref{count strings}) gives
$$v_{J_n} = \lim_{i \to \infty}\frac{\sum_{l \in \mathbb{A}}((\mathbf{A}^{i})_{J_n,l})}{\sum_{l \in \mathbb{A}}(\Lambda^i u_{l})}
= \lim_{i \to \infty}\frac{\sum_{l \in \mathbb{A}}((\mathbf{A}^{i})_{K_n,l})}{\sum_{l \in \mathbb{A}}(\Lambda^i u_{l})}= v_{K_n}.$$
By a similar argument using the equality of $J^{i-}$ and $K^{i-}$, we can deduce that $u_{J_1}=u_{K_1}$.  Finally, using the definition of the parry measure, 
$$\nu(cyl(j))=u_{J_1}v_{J_1}(v_{J_2}/(\Lambda v_{J_1}))...(v_{J_n}/(\Lambda v_{J_{n-1}}))=u_{J_1}v_{J_n}\Lambda^{n-1}$$
$$=u_{K_1}v_{K_n}\Lambda^{n-1} = u_{K_1}v_{K_1}(v_{K_2}/(\Lambda v_{K_1}))...(v_{K_n}/(\Lambda v_{K_{n-1}}))=\nu(cyl(k)) \blacksquare$$

This theorem allows us to realize the equivalence of measure between many cylinder sets once we make the following observation:

\begin{cor} \label{same measure}
If $j=J_1...J_n$ and $k=K_1...K_n$ are strings of lengthe $n$ such that $\Pi(J_1)=\Pi(K_1)$ and $\Pi(J_n)=\Pi(K_n)$, then $\nu(cyl(j))=\nu(cyl(k))$.
\end{cor}
 
Proof: Let $\mathbb{G}_2$ be the set $\{$A, A', A'', B, C, C', D, E$\}$ of vertices of $G_2$ excluding St and Fi.  Given $L_0 \in \mathbb{A}$, the number of strings of the form $L_0L_1$, such that $\Pi(L_1)=X$ for some fixed $X \in \mathbb{G}_2$, is simply the number on the arrow emanating from $\Pi(L_0)$ and terminating on $X$.  Similarly the number of strings of the form $L_{-1}L_0$, where $\Pi(L_{-1})=Y$ for some fixed $Y \in \mathbb{G}_2$, only depends on $\Pi(L_0)$ and $Y$.  By induction, we can see that the hypotheses of Theorem \ref{compare} on $J^{i+}, K^{i+}, J^{i-},$ and $K^{i-}$ are satisfied so we are done $\blacksquare$.

We will need one more lemma that states intuitively that "sharp turns" limit future possibilities.
\begin{lem} \label {contrast}
In the language of Theorem \ref {compare}, if $J^{i+} \ge K^{i+}$ for all $i$, and $\Pi{J_1}=\Pi{K_1}$, then $\nu(cyl(j)) \ge \nu(cyl(k))$.  In particular, if $J_1J_2$ and $K_1K_2$ are strings such that $\Pi{J_1}=\Pi{K_1} \in \{\text{A}, \text{A'}\}, \Pi(J_2) = \text{A},$ and $\Pi{K_2} \in \{\text{A'}, \text{A''}\}$, then $\nu(cyl(J_1J_2)) \ge \nu(cyl(K_1K_2))$.
\end{lem}

Proof: modifying the proof of Theorem \ref{compare}, (\ref{count strings}) becomes
$$\sum_{l \in \mathbb{A}}(\mathbf{A}^i)_{J_n,l} = J^{i+} \ge K^{i+} = \sum_{l \in \mathbb{A}}(\mathbf{A}^i)_{K_n,l},$$
so we have
$$v_{J_n} = \lim_{i \to \infty}\frac{\sum_{l \in \mathbb{A}}((\mathbf{A}^{i})_{J_n,l})}{\sum_{l \in \mathbb{A}}(\Lambda^i u_{l})} \ge \lim_{i \to \infty}\frac{\sum_{l \in \mathbb{A}}((\mathbf{A}^{i})_{K_n,l})}{\sum_{l \in \mathbb{A}}(\Lambda^i u_{l})}= v_{K_n}.$$ 
As in the proof of Corollary \ref{same measure}, $J^{i-} = K^{i-}$ for all $i$, so $u_{J_1} =  u_{K_1}$.  Thus 
$$\nu(cyl(j))=u_{J_1}v_{J_1}(v_{J_2}/(\Lambda v_{J_1}))...(v_{J_n}/(\Lambda v_{J_{n-1}}))=u_{J_1}v_{J_n}\Lambda^{n-1}$$
$$\ge u_{K_1}v_{K_n}\Lambda^{n-1} = u_{K_1}v_{K_1}(v_{K_2}/(\Lambda v_{K_1}))...(v_{K_n}/(\Lambda v_{K_{n-1}}))=\nu(cyl(k)).$$

For the special case, WLOG (by Corollary \ref{same measure}), we can let $J_2=$A$_a$ and $K_2\in\{$A'$^L_a$A''$^L_a\}$.  First consider the case where $K_2 = $A'$^L_a$.  For any string of the form $K_1K_2$A''$^L_bK_4...K_{i+2}$, $J_1J_2$A'$^L_bK_4...K_{i+2}$ is also a valid string.  Also, when $K_3 \ne $A''$^L_b$, if $K_1...K_{i+2}$ is a string, then so is $J_1J_2K_3...K_{i+2}$.  Since all of these strings beginning in $J_1J_2$ are different (whatever $K_3$ is), we know that there at least as many strings of length $i+2$ beginning in $J_1J_2$ as there are beginning in $K_1K_2$.

The case where $K_2 = $A''$^L_a$ is simpler.  For any string of the form $K_1...K_{i+2}$, $J_1J_2K_3...K_{i+1}$ is a different string.  In either case, $J^{i+} \ge K^{i+}$ so the result follows $\blacksquare$

\subsection{$H$ is not cohomologous to a constant} 
Once we have the following lemma about $H$, we will be in a position to prove Theorem \ref{variance} using a few known results.

\begin{lem} \label{cohomologous}
$H$ is not a cohomologous to a constant with respect to the operator induced by the shift map $\tau$.
\end{lem}

Proof: Suppose that $H$ is cohomologous to a constant i.e. suppose there exists a bounded function $H_0$ and constant $c$ such that $H(\mathbf{z})=c + H_0(\tau\mathbf{z})-H_0(\mathbf{z})$.  Then we have that
$$\sum_{i=0}^{n-1}H(\tau^i\mathbf{z})=\sum_{i=0}^{n-1}(c + H_0(\tau^{i +1}\mathbf{z})-H_0(\tau^i\mathbf{z}))=nc+H_0(\tau^n\mathbf{z})-H_0(\mathbf{z}).$$
Further, for any $n$ and any $\mathbf{x}$ and $\mathbf{y}$,

\begin{equation} \label{grows unbounded}
(\sum_{i=0}^{n-1}H(\tau^i\mathbf{x}))-(\sum_{i=0}^{n-1}H(\tau^i\mathbf{y}))=(H_0(\tau^n\mathbf{x})-H_0(\mathbf{x}))-(H_0(\tau^n\mathbf{y})-H_0(\mathbf{y}))
\end{equation}

is bounded.  Thus, in order to prove the lemma, it suffices to construct $\mathbf{x}$ and $\mathbf{y}$ where the quantity on the left grows unbounded with $n$.

Let $\mathbf{x} = ...$A$_a$A$_C|$A$_a$A$_C...$ (the vertical bar is a place-holder indicating the seperation between state $0$ and state $1$), representing a path that follows a geodic of the net $N$, cutting through the middle of successive copies of the fundamental domain.  Let $\mathbf{y} = ...$A'$^R_a$A'$^L_b|$A'$^R_a$A'$^L_b...$, representing a zig-zagging path in $N$ that alternates between sharp right turns and sharp left turns.  We will show that $H(\tau^i\mathbf{x}) - H(\tau^i\mathbf{y})$ is the same positive constant for all $i$.  Since $h_k \in \{0,1\}$, we can write $H(\mathbf{z})=\sum_{k=1}^{\infty}\nu\{\mathbf{Z}|h_k(\mathbf{z},\mathbf{Z})=1\}$.  For simplicity, we denote by $h_k\mathbf{z}$ the set $\{\mathbf{Z}|h_k(\mathbf{z},\mathbf{Z})=1\}$, so that $H(\mathbf{z})=\sum_{k=1}^{\infty}\nu(h_k\mathbf{z}).$

We wish to show that $H(\mathbf{x}) - H(\mathbf{y}) > 0$ so we will partition them as $H(\mathbf{x})=\sum_{j=1}^{12}p_j$ and  $H(\mathbf{y})=\sum_{j=1}^{11}q_j$ and observe (using Corollary \ref{same measure}) $p_j \ge q_j$ for $j \le10$ and $p_{11} > 0$.  Define the $p_j$'s and $q_j$'s as follows:

$$p_1=q_1=\nu(h_1\mathbf{x})=\nu(h_1\mathbf{y})= \nu (\bigcup_{\phi(w) \in \{\epsilon_a,\epsilon_A\}} cyl(w)),$$
\\
$$p_2=q_2=\nu((\bigcup_{w \in \{\text{A'}^R_B, \text{B}_{aB}, \text{A}_a, \text{B}_{da}, \text{A}_d, \text{B}_{Cd}\}} cyl(\text{D}_{Ba}w))$$
$$ \bigcup (\bigcup_{\phi(w) \in \{\alpha_b, \delta_{bA}, \alpha_A, \delta_{AD}, \alpha_D, \delta_{Dc}\}, w \ne \text{A''}^L_b} cyl(w\text{B}_{Ab}))),
$$
\\
$$ p_3=\nu((\bigcup_{w \in \{\text{A}_B, \text{B}_{aB}, \text{A}_a, \text{B}_{da}, \text{A'}^L_d\}} cyl(\text{D}_{cd}w))$$
$$ \bigcup (\bigcup_{\phi(w) \in \{\alpha_b, \delta_{bA}, \alpha_A, \delta_{AD}, \alpha_D\}, w \ne \text{A''}^R_D} cyl(w\text{B}_{DC})))
$$
$$=\nu((\bigcup_{w \in \{\text{A}_C, \text{B}_{DC}, \text{A}_D, \text{B}_{cD}, \text{A'}^L_c\}} cyl(\text{D}_{Ba}w))$$
$$ \bigcup (\bigcup_{\phi(w) \in \{\alpha_c, \delta_{cd}, \alpha_d, \delta_{dC}, \alpha_C\}, w \ne \text{A''}^R_C} cyl(w\text{B}_{Ab}))) = q_3$$
\\
$$p_4=\nu(cyl(\text{D}_{ab}\text{A'}^L_b))$$
$$=\sum_{j=0}^{\infty}(\nu(\bigcup_{w \ne \text{A'}^R_a}cyl(\text{D}_{ab}\text{A'}^L_b(\text{A'}^R_a\text{A'}^L_b)^jw))+\nu(\bigcup_{w \ne \text{A'}^L_b}cyl(\text{D}_{ab}(\text{A'}^L_b\text{A'}^R_a)^{j+1}w)))$$
$$\ge \sum_{j=0}^{\infty}\nu(\bigcup_{w \ne \text{A'}^R_a}cyl(\text{D}_{ab}\text{A'}^L_b(\text{A'}^R_a\text{A'}^L_b)^jw))=q_4$$
\\
$$p_5=\nu(\bigcup_{\phi(w)=\alpha_B}cyl(w\text{B}_{BA}))$$
$$=\sum_{j=0}^{\infty}(\nu(\bigcup_{\phi(w)=\zeta(\alpha_B\alpha_A)^j\alpha_B\beta_{BA}|\zeta \ne \alpha_A}cyl(w))+\nu(\bigcup_{\phi(w)=\zeta(\alpha_A\alpha_B)^{j+1}\beta_{BA}|\zeta \ne \alpha_B}cyl(w)))$$
$$\ge \sum_{j=0}^{\infty}\nu(\bigcup_{\phi(w)=\zeta(\alpha_B\alpha_A)^j\alpha_B\beta_{BA}|\zeta \ne \alpha_A}cyl(w))=q_5$$
\\
$$p_6=\nu(\bigcup_{\phi(w) \in \{\alpha_b,\delta_{bA},\alpha_A,\delta_{AD},\alpha_D\}=I_1}cyl(w\text{A}_b))$$
$$=\sum_{j=0}^{\infty}(\nu(\bigcup_{\phi(w)=\zeta\alpha_b(\alpha_a\alpha_b)^j\xi|\zeta \in I_1, \xi \ne \alpha_a}cyl(w))+\nu(\bigcup_{\phi(w)=\zeta(\alpha_b\alpha_a)^{j+1}\xi|\zeta \in I_1, \xi \ne \alpha_b}cyl(w)))$$
$$\ge \sum_{j=0}^{\infty}\nu(\bigcup_{\phi(w)=\zeta\alpha_b(\alpha_a\alpha_b)^j\xi|\zeta \in I_1, \xi \ne \alpha_a}cyl(w))=q_6$$
\\
$$p_7=\nu(\bigcup_{\phi(w) \in \{\alpha_B,\beta_{aB},\alpha_a,\beta_{da},\alpha_d\}=I_2, \phi(w')=\alpha_B} cyl(w'w))$$
$$=\sum_{j=0}^{\infty}(\nu(\bigcup_{\phi(w)=\zeta(\alpha_B\alpha_A)^j\alpha_B\xi|\zeta \ne \alpha_A, \xi \in I_2}cyl(w))+\nu(\bigcup_{\phi(w)=\zeta(\alpha_A\alpha_B)^{j+1}\xi|\zeta \ne \alpha_B, \xi \in I_2}cyl(w)))$$
$$\ge \sum_{j=0}^{\infty}\nu(\bigcup_{\phi(w)=\zeta(\alpha_B\alpha_A)^j\alpha_B\xi|\zeta \ne \alpha_A, \xi \in I_2}cyl(w))=q_7$$
\\
$$p_8=\nu(\bigcup_{\phi(w) \in \{\delta_{bA},\alpha_A,\delta_{AD},\alpha_D,\delta_{Dc}\},\phi(w')=\alpha_D}cyl(ww'))$$
$$=\nu(\bigcup_{\phi(w) \in \{\delta_{Dc},\alpha_c,\delta_{cd},\alpha_d,\delta_{dC}\}=I_3, \phi(w')=\alpha_b}cyl(ww'))$$
$$=\sum_{j=0}^{\infty}(\nu(\bigcup_{\phi(w)=\zeta\alpha_b(\alpha_a\alpha_b)^j\xi|\zeta \in I_3, \xi \ne \alpha_a}cyl(w))+\nu(\bigcup_{\phi(w)=\zeta(\alpha_b\alpha_a)^{j+1}\xi|\zeta \in I_3, \xi \ne \alpha_b}cyl(w)))$$
$$\ge \sum_{j=0}^{\infty}\nu(\bigcup_{\phi(w)=\zeta\alpha_b(\alpha_a\alpha_b)^j\xi|\zeta \in I_3, \xi \ne \alpha_a}cyl(w))=q_8$$
Here the second equality follows from Theorem \ref{compare}.
\\
$$p_9=\nu(\bigcup_{\phi(w) \in \{\beta_{Ba},\alpha_a,\beta_{da},\alpha_d,\beta_{Cd}\}, \phi(w')=\alpha_d} cyl(w'w))$$
$$=\nu(\bigcup_{\phi(w) \in \{\beta_{Cd},\alpha_C,\beta_{DC},\alpha_D,\beta_{cD}\}=I_4, \phi(w')=\alpha_B}cyl(w'w))$$
$$=\sum_{j=0}^{\infty}(\nu(\bigcup_{\phi(w)=\zeta(\alpha_B\alpha_A)^j\alpha_B\xi|\zeta \ne \alpha_A, \xi \in I_4}cyl(w))+\nu(\bigcup_{\phi(w)=\zeta(\alpha_A\alpha_B)^{j+1}\xi|\zeta \ne \alpha_B, \xi \in I_4}cyl(w)))$$
$$\ge \sum_{j=0}^{\infty}\nu(\bigcup_{\phi(w)=\zeta(\alpha_B\alpha_A)^j\alpha_B\xi|\zeta \ne \alpha_A, \xi \in I_4}cyl(w))=q_9$$
Here the second equality follows from Theorem \ref{compare}.
\\
$$p_{10}=\nu(\bigcup_{\phi(w) = \alpha_C, \phi(w') = \alpha_a} cyl(ww')) \ge \nu(\bigcup_{\phi(w) = \alpha_C, \phi(w') = \alpha_b} cyl(ww'))$$
$$=\sum_{j=0}^{\infty}(\nu(\bigcup_{\phi(w)=\alpha_C\alpha_b(\alpha_a\alpha_b)^j\xi| \xi \ne \alpha_a}cyl(w))+\nu(\bigcup_{\phi(w)=\alpha_C(\alpha_b\alpha_a)^{j+1}\xi| \xi \ne \alpha_b}cyl(w)))$$
$$\ge \sum_{j=0}^{\infty}\nu(\bigcup_{\phi(w)=\alpha_C\alpha_b(\alpha_a\alpha_b)^j\xi| \xi \ne \alpha_a}cyl(w))=q_{10}$$

Here the first inequality follows from four uses of Lemma \ref{contrast}.
\\
$$p_{11}=\nu(\bigcup_{\phi(w) = \alpha_A, \phi(w') = \alpha_c} cyl(ww')) \ge \nu(\bigcup_{\phi(w) = \alpha_B, \phi(w') = \alpha_c} cyl(ww'))$$
$$=\sum_{j=0}^{\infty}(\nu(\bigcup_{\phi(w)=\zeta(\alpha_B\alpha_A)^j\alpha_B\alpha_c|\zeta \ne \alpha_A}cyl(w))+\nu(\bigcup_{\phi(w)=\zeta(\alpha_A\alpha_B)^{j+1}\alpha_c|\zeta \ne \alpha_B}cyl(w)))$$
$$\ge \sum_{j=0}^{\infty}\nu(\bigcup_{\phi(w)=\zeta(\alpha_B\alpha_A)^j\alpha_B\alpha_c |\zeta \ne \alpha_A}cyl(w))=q_{11}$$
Again here the first inequality follows from four uses of Lemma \ref{contrast}.  Finally Let
$$p_{12} = H(\mathbf{x}) - \sum_{j=1}^{11}p_j > (\sum_{k=3}^{\infty}\nu(h_k\mathbf{x}))+\nu(cyl(\text{A}_c\text{A}_b)) =c_0$$
where $c_0$ is some positive constant.  Thus $H(\mathbf{x}) - H(\mathbf{y}) \ge c_0$.
By repeated uses of Theorem \ref{compare}, we can see that $H(\mathbf{x}) = H(\tau^i\mathbf{x})$ and $H(\mathbf{y}) = H(\tau^i\mathbf{y})$ for all positive integers $i$.  Thus the left-hand side of (\ref{grows unbounded})
$$(\sum_{i=0}^{n-1}H(\tau^i\mathbf{x}))-(\sum_{i=0}^{n-1}H(\tau^i\mathbf{y})) \ge n(c_0)$$ grows unbounded in $n$ yeilding the desired contradiction $\blacksquare$

We are now ready to prove Theorem \ref{variance}.  According to \cite{chas-lalley}, the varance of the distribution in Theorem \ref{Main} is 
$$\sigma^2 = \lim_{n \to \infty}(1/n)E(\sum_{i=1}^n\sum_{k=1}^{\infty}H_k(\tau^i\mathbf{X})-n\kappa)^2 = \lim_{n \to \infty}(1/n)E(\sum_{i=1}^nH(\tau^i\mathbf{X})-n\kappa)^2,$$
where
$$\kappa = \sum_{k=1}^{\infty}EH_k(\mathbf{X})=EH(\mathbf{X}).$$
Since $H$ does not depend on the past, this $\sigma^2$ is the same if calculated on the one-sided markov chain whose measure $\nu^{+}$ is determined by the Parry measure $\nu$.  $\nu^{+}$, being the unique measure of maximal entropy, is the equilibrium measure under the Reulle operator of the function $f=0$.  Under this measure, since $H$ isn't cohomologous to a constant by Lemma \ref{cohomologous}, Proposition 4.12 of Pollicott and Parry \cite{pp} gives $\sigma^2>0\blacksquare$

\end{document}